\documentclass[aap]{imsart}

\RequirePackage{amsthm,amsmath,amsfonts,amssymb}
\RequirePackage[authoryear]{natbib}
\RequirePackage{graphicx}

\usepackage[usenames,dvipsnames]{xcolor}

\usepackage{hyperref}
\hypersetup{pdfpagemode=FullScreen,  
	colorlinks=true,
	citecolor=Bittersweet,
	linkcolor=Bittersweet,
	urlcolor=Bittersweet
}

\startlocaldefs
\theoremstyle{plain}

\newtheorem{thm}{Theorem}[section]
\newtheorem{lem}{Lemma}[section]
\newtheorem{prop}{Proposition}[section]
\newtheorem{co}{Corollary}[section]

\newcommand{\cqfd}{\hfill $\square$}

\newcommand{\R}{\mathbb R}

\DeclareMathOperator*{\argmax}{\arg\max} 

\endlocaldefs

\begin{document}

\begin{frontmatter}
		\title{On the robustness of semi-discrete optimal transport}
		\runtitle{On the robustness of semi-discrete optimal transport}

\begin{aug}
			\author[A]{\fnms{Davy}~\snm{Paindaveine}\ead[label=e1]{Davy.Paindaveine@ulb.be}}
			\and
			\author[B]{\fnms{Riccardo}~\snm{Passeggeri}\ead[label=e2]{Riccardo.Passeggeri@imperial.ac.uk}}
			\address[A]{ECARES and Department of Mathematics,
				Universit\'{e} libre de Bruxelles, Brussels, Belgium\printead[presep={,\ }]{e1}}
			
			\address[B]{Department of Mathematics,
				Imperial College London, London, United Kingdom\printead[presep={,\ }]{e2}}
\end{aug}

\begin{abstract}
			We derive the breakdown point for solutions of semi-discrete optimal transport problems, which characterizes the robustness of the multivariate quantiles based on optimal transport proposed in \cite{GS}. We do so under very mild assumptions: the absolutely continuous reference measure is only assumed to have a support that is \textcolor{black}{convex}, whereas the target measure is a general discrete measure on a finite number, $n$ say, of atoms. The breakdown point depends on the target measure only through its probability weights (hence not on the location of the atoms) and involves the geometry of the reference measure through the \cite{Tuk1975} concept of  halfspace depth. Remarkably, depending on this geometry, the breakdown point of the optimal transport median can be strictly smaller than the breakdown point of the univariate median or the breakdown point of the spatial median, namely~$\lceil n/2\rceil /2$. In the context of robust location estimation, our results provide a subtle insight on how to perform multivariate trimming when constructing trimmed means based on optimal transport.
	\end{abstract}


\begin{keyword}
			\kwd{Breakdown point}
			\kwd{center-outward quantiles}
			\kwd{halfspace depth}
			\kwd{optimal transport}
			\kwd{robustness}
	\end{keyword}

\end{frontmatter}


	\section{Introduction}

	Dating back to \cite{Monge1781} and  \cite{Kanto1942}, the concept of optimal transport (OT) has found applications in very diverse fields, including economics (\citealp{Galichon2016}), 
	machine learning (\citealp{PeyrCutu2019}),
	signal processing (\citealp{Kolouri2017}), and computational biology (\citealp{Orlova16}), to mention only a few. Well-known monographs on the topic are, e.g., \cite{Villani2008}, \cite{Santa2015}, \cite{Panaretos2020}, and \cite{Ambrosio2021}.

	The literature dedicated to OT in probability and statistics has been exploding in the last decade. In particular, much effort has been dedicated to understanding the asymptotics of empirical OT costs; see, among many others, \cite{Boissard2014}, \cite{delBarrio2019}, \cite{WeedBach2019}, and \cite{HSM24}. Recently, \cite{Staudt2023} and  \cite{Manole2024a} tackled the unbounded setting, whereas  
	\cite{delBarrio2024} and \cite{Hundrieser2024b} studied in particular the semi-discrete case. 
	Of course, asymptotics for empirical OT maps themselves are also of primary interest; we refer, e.g., to \cite{Rigollet2021} and \cite{Sadhu2023} for convergence rates, and to \cite{Manole2024} for minimax rate optimality results; see also \cite{Bercu2021}, \cite{Goldfeld2024}, or \cite{Hundrieser2024} for results involving entropic OT.

	In another direction, OT has been the cornerstone to define modern concepts of multivariate quantiles and ranks; see \cite{Cheetal2017}, \cite{Hal21}, \cite{GS}, or, in the context of directional data, \cite{VerdeboutJRSSB} and \cite{BigotDirect24}. These new functionals opened the way to fully distribution-free rank tests in a very broad class of inference problems; see, among others, \cite{Bha2021}, 
	\cite{GS}, \cite{Shi2022}, \cite{ShiAoS22}, \cite{HallinJASA22,HallinBernou23}, \cite{DebSen2022}, \cite{HallinHlu2023}, and \cite{Shi24}.

	It has been recognized, however, that OT is sensitive to outliers. As a consequence, several robustifications of OT have been proposed in the literature; see, e.g., 
	\cite{Alvarez2008},
	\cite{Balayi2020}, 
	\cite{Le2021}, 
	%
	and
	\cite{Nietert2022}. In line with this, the minimum Kantorovich estimation procedure from \cite{Bassetti2006}, which is a parametric estimation method based on OT, was robustified in \cite{Balayi2020}; see also \cite{Mukherjee2021}. Recently, \cite{Ronch2023} commented on the (lack of) robustness of OT but also made the point that this had not yet been properly investigated/quantified in the literature. In particular, the author points out that \emph{``concepts of local stability such as the influence function and of global reliability such as the breakdown point still have to be developed''.}

    
\textcolor{black}{In the direction of influence functions, Gateaux derivatives of the OT map (involving a smooth perturbation measure rather than a Dirac one as in standard influence functions) were actually obtained in \cite{Loeper2005}, 
    \cite{ManoleBNW2024}, and \cite{GonSheng2024}. We also refer to \cite{Sadhu2023} for the semi-discrete setting, and to \cite{Gonz2022} and \cite{Goldfeld2024} for regularized OT. In contrast, no results are available for global robustness of OT in terms of \emph{breakdown point}; see \cite{Hametal1986} and \cite{Dono83}. This provides a natural motivation for the present work, that derives the breakdown point of \emph{semi-discrete}\footnote{\textcolor{black}{In an interesting independent work, \cite{Ave2024} actually tackled the same problem when the reference and target measures are rather both discrete or both  continuous; see Section~\ref{secPersp} below.}} OT}; see Section~\ref{secSemiDiscrete} for a precise definition. As a corollary, this will in particular provide the breakdown point of OT-based multivariate quantiles; we refer to \cite{AIHP} for a recent breakdown point analysis of a competing concept of multivariate quantiles, namely the \emph{geometric} or \emph{spatial} quantiles from \cite{DuKol92} and \cite{Cha1996}; see also \cite{Kol1997}.

	To some extent, the problem we consider in this paper is related to the stability analysis of solutions of OT problems, where the objective is to study how the solutions are affected when one slightly changes the target measure. A recent paper in this line of research is \cite{Bansil}, in which the authors fix the location of the atoms and investigate the behavior of the solution as the weights change, focusing on functional estimates (e.g.~$L_2$ or uniform). However, we cannot use their results as functional estimates do not translate into breakdown point results.


	\subsection{Semi-discrete OT} 
	\label{secSemiDiscrete}

	Our starting point is Monge's problem under the $L_2$ loss: given a probability measure $\mu$ on $\mathcal{S}\subset\mathbb{R}^d$ and a probability measure $\nu$ on $\mathcal{X}\subset\mathbb{R}^d$, one needs to find the measurable map $T$ from $\mathcal{S}$ to $\mathcal{X}$ solving
	\begin{equation}
		\label{Monge}
		\inf\limits_{T}\int\|u-T(u)\|^2 d\mu(u)\quad\text{subject to $T_{\#}\mu=\nu$}
		,
	\end{equation}
	where the pushforward measure~$T_{\#}\mu$ of $\mu$ by $T$ is the one that associates to any borel set~$B\subset\mathcal{X}$ the measure~$T_{\#}\mu(B)=\mu(T^{-1}(B))$. As usual, we call $\mu$ the \emph{reference measure} and $\nu$ the \emph{target measure}. In this framework, a result of paramount importance is the Brenier--McCann theorem; see \cite{Brenier} and \cite{McCann}. We state it here in its version given in \cite{GS}\textcolor{black}{, with the only difference that we will denote~$Q_\nu$ (rather than~$Q$) the optimal transport map.\footnote{\textcolor{black}{The reason is that our focus in the present work will be on the sensitivity of the optimal transport map with respect to the target measure~$\nu$, whereas the reference measure~$\mu$ will be fixed.}}}

	\begin{thm}[Brenier--McCann theorem]
		\label{TheorBM}
		Let $\mu$ and $\nu$ be Borel probability measures on $\mathbb{R}^d$. Suppose further that $\mu$ is absolutely continuous\footnote{Throughout, absolute continuity is with respect to the Lebesgue measure on~$\R^d$.}. Then there exists a convex
		function $\psi:\mathbb{R}^d\to\mathbb{R}\cup\{\infty\}$ whose gradient $Q_\nu = \nabla\psi : \mathbb{R}^d\to\mathbb{R}^d$ pushes $\mu$ forward to $\nu$.
		In fact, there exists only one such $Q_\nu$ that arises as the gradient of a convex function, that is, $Q_\nu$ is unique $\mu$-almost everywhere. Moreover, if $\mu$ and $\nu$ have finite second moments, then $Q_\nu$ uniquely minimizes Monge's problem~(\ref{Monge}).
	\end{thm}

	Throughout, we will assume that $\mu$ is an absolutely continuous probability measure with \textcolor{black}{convex} support\footnote{In this work, we define the support of a probability measure~$\mu$ on~$\R^d$ as the smallest closed set with $\mu$-measure one.} $\mathcal{S}$ and that~$\nu$ is a finitely discrete probability measure with atoms~$x_1,\ldots,x_n\in\mathbb{R}^d$ and weights~$\lambda_1,\ldots,\lambda_n$, that is,
	$$
	\nu
	=
	\sum_{i=1}^{n}\lambda_i\delta_{x_i}
	;
	$$ 
	as usual, $\delta_x$ here stands for the Dirac probability measure at~$x$. \textcolor{black}{We stress that we do not assume that the support of~$\mu$ is compact, which allows one to adopt, e.g., the $d$-variate standard normal distribution as a reference distribution---as was done in \cite{Tengyao2025} when conducting multiple-output composite quantile regression. Of course, as soon as~$\mu$ has finite second moments, then both~$\mu$ and~$\nu$ have, so that} Theorem~\ref{TheorBM} states that the unique gradient of a convex function that pushes~$\mu$ forward to~$\nu$ uniquely minimizes Monge's problem~(\ref{Monge}). A particular case of interest will be the case of \emph{empirical target measures}, that is obtained with~$\lambda_1=\ldots=\lambda_n=1/n$.   
	
	Since the reference measure is absolutely continuous while the target measure is discrete, this is usually referred to as the \emph{semi-discrete} OT framework. To describe the corresponding solution~$Q_\nu$ to Monge’s problem in this setup, we need to introduce the following concepts.
	Let $\mathcal{X}=\{x_1,\ldots,x_n\}$ be a finite set of points in $\mathbb{R}^d$, and $w=(w_1,\ldots,w_n)\in\mathbb{R}^d$ be a given vector (usually called \emph{weight vector} even though $w_1,\ldots,w_n$ may be negative). Then, for any~$i=1,\ldots,n$, the \emph{power cell} of~$x_i$ with respect to~$\mathcal{S}$ is defined as
	$$
	{\rm Lag}^{w}_\mathcal{X}(i)
	=
	\big\{
	u \in \mathcal{S} : 
	\| u-x_i\|^2-w_i\leq \| u-x_j\|^2-w_j
	\textcolor{black}{\textrm{ for } j\in\{1,\ldots,n\}}
	\big\}
\textcolor{black}{	;}
	$$
\textcolor{black}{these cells are often referred to as \emph{Laguerre cells}, which explains the notation.} The collection of power cells~${\rm Lag}^{w}_\mathcal{X}(i)$, $i=1,\ldots,n$, provides a decomposition of~$\mathcal{S}$ that is called the \emph{power diagram} or \emph{weighted Voronoi diagram} of $(\mathcal{X}, w)$ with respect to~$\mathcal{S}$; we refer to \cite{AHA87} for detailed properties of such power diagrams. The corresponding \emph{power map} $T^{w}_{\mathcal{X}}:\mathcal{S}\to\mathcal{X}$ is then such that $T^{w}_{\mathcal{X}}(u)=x_i$ if $u\in {\rm Lag}^{w}_\mathcal{X}(i)$. This map is well-defined except on the boundary of the power cells, hence it is well-defined $\mu$-almost everywhere. Now, a weight vector $w$ is said to be \textit{adapted} to the couple of measures $(\mu, \nu)$ if 
	\begin{equation*}
		\lambda_i
		=
		\mu({\rm Lag}^{w}_\mathcal{X}(i))
		\
		\textrm{ for~}i=1,\ldots,n
		.
	\end{equation*}
	As shown in Theorem~2.5 from~\cite{Lev}, such an adapted weight vector always exists. Crucially,  Theorem 2 in \cite{Mer}, which is based on Section 5 of \cite{AHA98} (see also \citealp{GS}), implies that if~$w$ is an adapted weight vector, then 
	$$
	Q_\nu
	=
	T^{w}_{\mathcal{X}}
	$$
	provides the solution to Monge’s problem (\ref{Monge}) in the considered semi-discrete framework. \textcolor{black}{Note that the adapted weight vector is unique up to an additive constant (see Section~3.4 of \citealp{Mer}) and is the optimal transport potential, namely it is the solution of the Kantorovich dual program, as noted in Section 2.5 of~\cite{Lev}}.

	\subsection{Breakdown point for semi-discrete OT}

	We now define the breakdown point of the solution to Monge's problem in the semi-discrete case, namely the  breakdown point of~$Q_\nu$. We will actually consider the breakdown point of~$Q_\nu(u)$ for any~$u\in\mathcal{S}$. To do so, denote as~$\mathcal{B}_\delta(u)$ the open ball with center~$u$ and radius~$\delta(>0)$. We then define the breakdown point of $Q_\nu(u)$ as
	\begin{eqnarray}
		\lefteqn{
			\hspace{1mm}
		{\rm BDP}(Q_\nu(u))
			=
			\min
			\Bigg\{
			\sum_{i\in I}\lambda_i: I\subset\{1,\ldots,n\}
			\ \ \mathrm{ such\ that } 
		}
		\nonumber
		\\[2mm]
		& & 
		\hspace{41mm}
		\sup
		\int_{\mathcal{B}_\delta(u)}
		\|Q_\nu(x)-Q_{\tilde{\nu}_I}(x)\| \, d\mu(x)
		=
		\infty
		\ \ \
		\forall \delta>0
		\Bigg\}
		,
		\label{defBDP}
	\end{eqnarray}
	where the supremum is taken over 
	target measures
	$$
	\tilde{\nu}_I=\sum_{j=1}^m\tilde{\lambda}_j\delta_{\tilde{x}_j}
	,
	$$
	where the atoms $\tilde{x}_1,\ldots,\tilde{x}_m\in\mathbb{R}^d$ include~$x_i$, $i\in I^c:=\{1,\ldots,n\}\setminus I$, and where the weights $\tilde{\lambda}_1,\ldots,\tilde{\lambda}_m$ are so that~$\tilde{\nu}_I(\{x_i\})=\lambda_i$ for any~$i\in I^c$. The observations~$x_i$, $i\in I^c$, are thus those that are not contaminated and
	$$
	\sum_{i\in I}\lambda_i
	$$
	is the total mass of contamination. Note that our definition of breakdown point is more general than the usual one from \cite{Hampel68,Hampel71}. In particular, we allow the contamination to affect not only the locations of certain atoms but also their weights. Also, we allow the contamination to modify the numbers of atoms considered; this may be of interest when considering, e.g., symmetric minimal contaminations as in Theorem~3.3 of \cite{AIHP}.
	The reason why we consider 
		$$
        \sup
		\int_{\mathcal{B}_\delta(u)}
		\|Q_\nu(x)-Q_{\tilde{\nu}_I}(x)\| \, d\mu(x)
		=
		\infty
		\quad
		\forall	 \delta>0
    $$
	in~(\ref{defBDP}), instead of just $\sup\|Q_\nu(u)-Q_{\tilde{\nu}_I}(u)\|=\infty$, is that $Q_\nu(u)$ and $Q_{\tilde{\nu}_I}(u)$ are defined only $\mu$-almost everywhere.  \textcolor{black}{Finally, note that the breakdown point~${\rm BDP}(Q_\nu(u))$ is local in~$u$ since it would be perfectly equivalently to require that the supremum in~(\ref{defBDP}) is infinite only for any~$\delta>0$ small enough}\textcolor{black}{---in the sense that there exists~$\delta_0>0$ (that may depend on~$I$ and on~$u$) such that the supremum in~(2) is infinite for any~$\delta\in(0,\delta_0)$.}


	\subsection{Outline 
		and notation}

	The outline of the paper is as follows. In Section~\ref{secMainResults}, we state our main results, and in particular Theorem~\ref{TheorBDP}, that provides a sharp expression for the breakdown point of~$Q_\nu$ at any point~$u$ in~$\mathcal{S}$. We discuss the breakdown point of the OT median and consider the case of OT quantiles for several classical reference measures~$\mu$.  We also comment on the relevance of our results in the construction of OT trimmed means. In Section~\ref{secProofUpperBound} (resp., Section~\ref{secProofLowerBound}), we establish the upper bound (resp., lower bound) for the breakdown point, which together prove Theorem~\ref{TheorBDP}. In Section~\ref{secPersp}, we provide \textcolor{black}{final comments} and briefly discuss perspectives for future research. Finally, the statements and proofs of some auxiliary results are provided in an appendix.
	%
	For the sake of convenience, we introduce here some notation we will use in the paper. The unit sphere in~$\R^d$ will be denoted as~$\mathcal{S}^{d-1}=\{x\in\R^d:\|x\|=1\}$. For a subset~$A$ of~$\R^d$, we will write~$A^c$, ${\rm int}(A)$, ${\rm cl}(A)$, and~$\partial A$ for the complement, the interior, the closure, and the boundary of~$A$, respectively. The  distance between~$x\in\R^d$ and~$A$ is defined as~$d(x,A):=\inf \{\|x-y\|:y\in A\}$. \textcolor{black}{For any~$x,y\in\R^d$, we will denote the inner product~$\sum_{i=1}^d x_i y_i$ as~$\langle x,y\rangle$.} The notation~$\mathcal{B}_r$ will denote the open ball with radius~$r$ centered at the origin of~$\R^d$. The Lebesgue measure on~$\R^d$ will be denoted as~$\mathcal{L}$. In the sequel, $\lceil \cdot\rceil$ will denote a slightly modified ceiling function: if~$t\neq 0$, then we define~$\lceil t\rceil$ as the usual ceiling of~$t$ (that is, as the smallest integer that is larger than or equal to~$t$), but we let~$\lceil 0\rceil=1$. Finally, a sum of an empty set of terms is defined as zero.


	\section{Main results}
	\label{secMainResults}

\textcolor{black}{Semi-discrete OT enjoys the natural translation-equivariance property described in Theorem~\ref{TheorTrivialUpperBound}(i) below.}
	By adapting the argument used in the proof of Theorem~2.1 from \cite{LopRou1991}, this equivariance property allows us to derive an upper bound for the breakdown point of~$Q_\nu(u)$. More precisely, we have the following result.

	\begin{thm}
		\label{TheorTrivialUpperBound}
		Let $\mu$ be an absolutely continuous probability measure over~$\R^d$ with
         \textcolor{black}{support}~$\mathcal{S}$ and let~$\nu$ be the discrete probability measure with atoms~$x_1,\ldots,x_n\in\mathbb{R}^d$ and weights~$\lambda_1,\ldots,\lambda_n$. Fix~$u\in \mathcal{S}$.
		Then, 
		\begin{enumerate}
		\item[(i)]	
$Q_\nu(u)$ is translation-equivariant in the sense that if
		$$
		\nu_t
		=
		\sum_{i=1}^{n}\lambda_i\delta_{x_i+t}
		$$
		for some $d$-vector~$t$, then 
		$
		Q_{\nu_t}(u)
		=
		Q_\nu(u)+t
		$.
				\item[(ii)]	
Consequently, the breakdown point of~$Q_\nu(u)$ satisfies
		$$
		{\rm BDP}(Q_\nu(u))
		\leq
		\min
		\Bigg\{
		\sum_{i\in I}\lambda_i: I\subset\{1,\ldots,n\}
		\ \, \mathrm{ such\ that }\, \
		\sum_{i\in I}\lambda_i\geq \frac{1}{2} 
		\Bigg\}
		.
		$$
			\item[(iii)]	
 For empirical target measures~$\nu=\nu_n$ say, this yields ${\rm BDP}(Q_{\nu_n}(u))
		\leq
		\left\lceil\frac{n}{2}\right\rceil / n$, so that 
		$$
		\limsup_{n\to\infty}
		{\rm BDP}(Q_{\nu_n}(u))
		\leq
		\frac{1}{2}
\cdot
$$
		\end{enumerate}
	\end{thm}

	\noindent
	{\sc Proof of Theorem~\ref{TheorTrivialUpperBound}.}
	(i) \textcolor{black}{The result is stated in Lemma~A.7 from \cite{GSsupp} and, as pointed out by one of the Reviewers, it follows directly from the fact that shifts of maximal monotone operators are maximal monotone.}
	%
	(ii)
	Fix~$I\subset\{1,\ldots,n\}$ with
	$$
	\sum_{i\in I}\lambda_i\geq \frac{1}{2} 
	,
	$$
	so that the total probability mass of the uncontaminated atoms satisfies
	$$
	\sum_{i\in I^c}
	\lambda_i
	\leq \frac{1}{2} 
	\cdot
	$$
	For any $d$-vector~$t$, consider then the contaminated target measure
	$$
	\tilde{\nu}(t)
	=
	\sum_{i\in I^c}
	\lambda_i
	\delta_{x_i}
	+
	\sum_{i\in I^c}
	\lambda_i
	\delta_{x_i+t}
	+
	\bigg(
	1
	-
	2\sum_{i\in I^c}
	\lambda_i
	\bigg)
	\delta_{t/2}
	.
	$$
	Note that, for sufficiently large~$\|t\|$, the contaminated measure~$\tilde{\nu}(t)$ attributes measure~$\lambda_i$ to~$\{x_i\}$ for any~$i\in I^c$, so that~$\tilde{\nu}(t)$ is among the contaminated measures on which the supremum is considered in~(\ref{defBDP}). Now, since
	$$
	\tilde{\nu}(-t)
	=
	\sum_{i\in I^c}
	\lambda_i
	\delta_{x_i-t}
	+
	\sum_{i\in I^c}
	\lambda_i
	\delta_{x_i}
	+
	\bigg(
	1
	-
	2\sum_{i\in I^c}
	\lambda_i
	\bigg)
	\delta_{-t/2}
	,
	$$
	translation-equivariance implies that, for any~$\delta>0$, 
	$$
	\int_{\mathcal{B}_\delta(u)}
	\|Q_{\tilde{\nu}(t)}(x)-Q_{\tilde{\nu}(-t)}(x)\| 
	\, d\mu(x)
	=
\textcolor{black}{
	\int_{\mathcal{B}_\delta(u)}
	\|t\|
	\, d\mu(x)
	=
	\|t\|
}
	\mu(\mathcal{B}_\delta(u))
	$$
	diverges to infinity as~$\|t\|$ does (\textcolor{black}{we must have~$\mu(\mathcal{B}_\delta(u))>0$ because if we would have~$\mu(\mathcal{B}_\delta(u))=0$, then}~$\mathcal{S}\setminus \mathcal{B}_\delta(u)$ would be a closed subset of~$\mathcal{S}$ with~$\mu$-measure one, which would contradict the fact that~$\mathcal{S}$ is the support of~$\mu$). However, \textcolor{black}{if by contradiction we would} have~${\rm BDP}(Q_\nu(u))>\sum_{i\in I}\lambda_i$, then, for sufficiently large~$\|t\|$, we would have
	\begin{eqnarray*}
		\lefteqn{
			\int_{\mathcal{B}_\delta(u)}
			\|Q_{\tilde{\nu}(t)}(x)-Q_{\tilde{\nu}(-t)}(x)\| 
			\, d\mu(x)
		}
		\\[2mm]
		&  &
		\hspace{13mm}
		\leq 
		\int_{\mathcal{B}_\delta(u)}
		\|Q_\nu(x)-Q_{\tilde{\nu}(t)}(x)\| 
		\, d\mu(x)
		+
		\int_{\mathcal{B}_\delta(u)}
		\|Q_\nu(x)-Q_{\tilde{\nu}(-t)}(x)\| 
		\, d\mu(x)
		\\[2mm]
		&  &
		\hspace{13mm}
		\leq 
		2
		\sup 
		\int_{\mathcal{B}_\delta(u)}
		\|Q_\nu(x)-Q_{\tilde{\nu}_I}(x)\| 
		\, d\mu(x)
		\\[2mm]
		&  &
		 \hspace{13mm}
		<
		\infty
	\end{eqnarray*}
	\textcolor{black}{for some~$\delta>0$} (here, we used the fact that, for sufficiently large~$\|t\|$, $\tilde{\nu}(-t)$ is among the contaminated measures on which the supremum is considered in~(\ref{defBDP}), too). Since this is a contradiction, we must have~${\rm BDP}(Q_\nu(u))\leq \sum_{i\in I}\lambda_i$ for any~$I\subset\{1,\ldots,n\}$ with
	$
	\sum_{i\in I}\lambda_i\geq 1/2
	$,
	which establishes the result. Of course, (iii) is a direct corollary of~(ii), so that the result is proved. 
	\cqfd
	\vspace{3mm}


	Part~(ii) of this result might suggest that, as it is the case for the usual quantiles in dimension~$d=1$ and for the spatial quantiles in arbitrary dimension~$d$ (see \citealp{AIHP}), the breakdown point of the quantile~$Q_\nu(u)$ has a maximum value (in~$u$) that is equal to~$\lceil n/2\rceil/n$. Interestingly, it turns out that this is the case for some reference measures~$\mu$ only. This will be one of the many consequences of the following result,  which provides an explicit expression for the breakdown point of~$Q_\nu(u)$. 
	\textcolor{black}{Interestingly, this expression depends on~$u$ only through its \cite{Tuk1975} \emph{halfspace depth} with respect to~$\mu$, that is defined as
		$H\!D(u,\mu) 
		=
		\inf
		\{
		\mu(H): H\ \mathrm{a\ closed\ halfspace\ containing\ } u 
		\}  
		$}.

	\begin{thm}
		\label{TheorBDP}
		Let $\mu$ be an absolutely continuous probability measure over~$\R^d$ with \textcolor{black}{convex} support $\mathcal{S}$ and let~$\nu$ be the discrete probability measure with atoms~$x_1,\ldots,x_n\in\mathbb{R}^d$ and weights~$\lambda_1,\ldots,\lambda_n$. Fix~$u\in \mathcal{S}$. Then, 
		\begin{enumerate}
		\item[(i)]	
		the breakdown point of~$Q_\nu(u)$ is
		$$
		{\rm BDP}(Q_\nu(u))
		=
		\min
		\Bigg\{
		\sum_{i\in I}\lambda_i: \emptyset \neq I\subset\{1,\ldots,n\}
		\ \, \mathrm{ such\ that }\, \
		\sum_{i\in I}\lambda_i\geq H\!D(u,\mu) 
		\Bigg\}
		.
		$$ 
				\item[(ii)]	
 For empirical measures~$\nu$, this yields 
		\begin{equation*}
			{\rm BDP}(Q_\nu(u))
			=
			\left\lceil n\, H\!D(u,\mu) \right\rceil / n
		\end{equation*}
		$($recall that we defined~$\lceil 0\rceil=1)$. In particular, ${\rm BDP}(Q_\nu(u))\to H\!D(u,\mu)$ as $n\to\infty$.
		\end{enumerate}
	\end{thm}


	%
	%
	%
	%

	A direct consequence of Theorem~\ref{TheorBDP} is that the breakdown point of~$Q_\nu(u)$ is maximized at the \emph{Tukey median of~$\mu$}, that is defined as 
	$$
	u_*
	=
	\argmax_{u\in \mathcal{S}}
	H\!D(u,\mu)
	;
	$$
	existence of a maximizer holds for any probability measure~$\mu$ (see, e.g., Proposition~7 in \citealp{RouRut1999}), and its uniqueness is guaranteed if any neighborhood of~$u_*$ has a positive $\mu$-measure (which is the case under the assumptions we consider on~$\mu$). The corresponding semi-discrete OT quantile~$Q_\nu(u_*)$ is then a natural definition for the \emph{OT median} of~$\nu$. The following result concerns the breakdown point of this median.

	\begin{co}
		\label{CorollaryBDPdepth}
		Let $\mu$ be an absolutely continuous probability measure over~$\R^d$ with \textcolor{black}{convex} support $\mathcal{S}$ and let~$\nu$ be the discrete probability measure with atoms~$x_1,\ldots,x_n\in\mathbb{R}^d$ and weights~$\lambda_1,\ldots,\lambda_n$. Then, 
		\begin{enumerate}
				\item[(i)]	
 the maximum breakdown point of~$Q_\nu(u)$ over~$\mathcal{S}$ is
		$$
		\max_{u\in \mathcal{S}}
		{\rm BDP}(Q_\nu(u))
		=
		{\rm BDP}(Q_\nu(u_*))
		,
		$$
		and this maximal breakdown point satisfies
		$$
		{\rm BDP}(Q_\nu(u_*))
		\geq 
		\min
		\Bigg\{
		\sum_{i\in I}\lambda_i: I\subset\{1,\ldots,n\}
		\ \, \mathrm{ such\ that }\, \
		\sum_{i\in I}\lambda_i\geq \frac{1}{d+1} 
		\Bigg\}
		.
		$$
						\item[(ii)]	
If~$\mu$ is angularly symmetric about~$\theta$ (in the sense that~$\mu(\theta-B)=\mu(\theta+B)$ for any Borel cone~$B$, that is, for any Borel set~$B$ such that~$rB=B$ for any~$r>0$), then~$u_*=\theta$ and
		$$
		{\rm BDP}(Q_\nu(u_*))
		=
		\min
		\Bigg\{
		\sum_{i\in I}\lambda_i: I\subset\{1,\ldots,n\}
		\ \, \mathrm{ such\ that }\, \
		\sum_{i\in I}\lambda_i\geq \frac{1}{2} 
		\Bigg\}
		;
		$$
		in particular, for empirical target measures~$\nu=\nu_n$, we then have
		$
		{\rm BDP}(Q_{\nu_n}(u_*))
		=
		\left\lceil \frac{n}{2} \right\rceil / n
		$, so that ${\rm BDP}(Q_{\nu_n}(u_*))
		%
		%
		%
		%
		%
		%
		\to 1/2$ as $n\to\infty$.
		\end{enumerate}
	\end{co}

	\noindent
	{\sc Proof of Corollary~\ref{CorollaryBDPdepth}.}
	In view of Theorem~\ref{TheorBDP}, Part~(i) of the result is a direct consequence of \textcolor{black}{Lemma~6.3} from \cite{DonGas1992} \textcolor{black}{(or, alternatively, of Proposition~9 from \citealp{RouRut1999})}. 
	\color{black}
	We thus turn to Part~(ii). First, since~$\mu$ is absolutely continuous and angularly symmetric about~$\theta$, Theorem~1 in \cite{RouStr2004} implies that~$H\!D(\theta,\mu)=1/2$, so that Corollary~1 in the same paper then entails that~$H\!D(\theta,\mu)$ maximizes~$H\!D(\cdot,\mu)$ over~$\R^d$ (equivalently, over~$\mathcal{S}$, since~$H\!D(x,\mu)=0$ for any~$x\notin \mathcal{S}$). It remains to show that~$u_*=\theta$ is the only maximizer of~$H\!D(\cdot,\mu)$ over~$\R^d$. \emph{Ad absurdum}, assume that~$\xi\in\R^d\setminus\{\theta\}$ satisfies~$H\!D(\xi,\mu)=1/2$. Then, since~$\mu$ is absolutely continuous, Theorem~2 in \cite{RouStr2004} implies that~$\mu$ is angularly symmetric about~$\xi$, hence is angularly symmetric about both~$\theta$ and~$\xi$. Lemma~2.3 in \cite{ZuoSer2000D} then entails that~$\mu$ is concentrated on a line of~$\R^d$, which contradicts absolute continuity of~$\mu$.   
	\color{black}
	\cqfd
	\vspace{3mm}

	Interestingly, it is \emph{only} for angularly symmetric probability measures~$\mu$ that we would have ${\rm BDP}(Q_{\nu_n}(u_*))\to 1/2$ as $n\to\infty$ for empirical target measures~$\nu$; this is a corollary of Theorem~2 in \cite{RouStr2004}. For any~$\eta>0$, it is easy to construct a reference measure~$\mu$ satisfying the assumptions of Theorem~\ref{TheorBDP} for which $\lim_{n\to\infty} {\rm BDP}(Q_{\nu_n}(u_*))\in [1/(d+1),1/(d+1)+\eta]$ for empirical target measures~$\nu$. In other words, the asymptotic breakdown point of the OT median can be arbitrarily close to~$1/(d+1)$.   
	
	Now, while Corollary~\ref{CorollaryBDPdepth} focuses on the breakdown point of the innermost OT quantile, it is natural to consider other quantiles, too, which we now do. For the sake of clarity, we will restrict here to empirical target measures~$\nu_n$, but of course the result could be stated for arbitrary finitely discrete target measures. The reference measure~$\mu$ mainly considered in \cite{GS} is the uniform distribution on the hypercube~$[0,1]^d$. For~$d=2$, it directly follows from Section~5.4 in \cite{RouRut1999} that, 	for any~$u\in(u_1,u_2)\in[0,1]^2$,  
	\begin{equation*}
		{\rm BDP}(Q_{\nu_n}(u))
		=
		\big\lceil 2n\, \min(u_1,1-u_1)\min(u_2,1-u_2) \big\rceil / n
		.
	\end{equation*}
	As mentioned in Example~1 of \cite{NagySS}, the corresponding expression for~$d>2$, that is much more involved, can be obtained from Lemma~1.3 (and its proof) in \cite{Schu1991}. Other reference measures, which allow for a natural, directional, indexing of OT quantiles, are measures~$\mu$ that are supported on~$\mathcal{S}={\rm cl}(\mathcal{B}_1)$ and are orthogonal-invariant (in the sense that $\mu(OB)=\mu(B)$ for any Borel set~$B$ in~$\R^d$ and any $d\times d$ orthogonal matrix~$O$). An example is the \emph{spherical uniform measure} on~${\rm cl}(\mathcal{B}_1)$, that is the orthogonal-invariant measure~$\mu$ for which~$\mu(\mathcal{B}_r)=r$ for any~$r\in[0,1]$; see \cite{Cheetal2017} and~\cite{Hal21}. Of course, another example is the uniform measure on~${\rm cl}(\mathcal{B}_1)$. We have the following result.

	\begin{co}
		\label{CorollarySpherQuantiles}
		Assume that the reference measure~$\mu$ is orthogonal-invariant and has support~$\mathcal{S}={\rm cl}(\mathcal{B}_1)$. Let~$\nu_n$ be an empirical target measure. Then,
\begin{enumerate}
						\item[(i)]	
 for any~$\alpha\in[0,1]$ and any~$v\in\mathcal{S}^{d-1}$,
		\begin{equation*}
			{\rm BDP}(Q_{\nu_n}(\alpha v))
			=
			\left\lceil n \mu( \mathcal{B}_1 \cap \{u\in\R^d: \alpha \leq u_1 \leq 1 \}) \right\rceil / n
			.
		\end{equation*}
		%
								\item[(ii)]	
 If~$\mu$ is the spherical uniform measure on~${\rm cl}(\mathcal{B}_1)$, then
		\begin{equation}
			\label{BDPspherunif}
			{\rm BDP}(Q_{\nu_n}(\alpha v))
			=
			\left\lceil  
			\frac{n\Gamma(\frac{d}{2})}{\sqrt{\pi}\Gamma(\frac{d-1}{2})}
			\int_\alpha^1
			\int_{x}^1
			r^{-(d-2)} 
			(r^2-x^2)^{(d-3)/2}
			\,
			dr
			\,
			dx
			\right\rceil \big/ n
			.
		\end{equation}
								\item[(iii)]	 
								If~$\mu$ is the uniform measure on~${\rm cl}(\mathcal{B}_1)$, then
		\begin{equation}
			\label{BDPunif}
			{\rm BDP}(Q_{\nu_n}(\alpha v))
			=
			\left\lceil  
			\frac{n\Gamma(\frac{d+2}{2})}{\sqrt{\pi}\Gamma(\frac{d+1}{2})}
			\int_\alpha^1
			(1-x^2)^{(d-1)/2}
			\,
			dx
			\right\rceil \big/ n
            ,
		\end{equation}
		where~$\Gamma$ is the Euler Gamma function.  
\end{enumerate}		
	\end{co}

	\noindent
	{\sc Proof of Corollary~\ref{CorollarySpherQuantiles}.}
	(i) Since~$\mu$ is orthogonal-invariant, $H\!D(u,\mu)=\mu(\{u\in\R^d: u_1\leq -\|u\|\})=\mu(\{u\in\R^d: u_1\geq \|u\|\})$; see, e.g., Example~2 in \cite{NagySS}. Since~$\mu$ has support~${\rm cl}(\mathcal{B}_1)$, the result then follows from Part~(ii) of Theorem~\ref{TheorBDP}. (ii) If~$X=(X_1,\ldots,X_d)$ follows the spherical uniform distribution~$\mu$, then~$\|X\|$ is uniformly distributed on~$[0,1]$, so that~$X_1$ admits the density
	$$
	f_1(x)
	=
	\bigg(
	\frac{\Gamma(\frac{d}{2})}{\sqrt{\pi}\Gamma(\frac{d-1}{2})}
	\int_{|x|}^1
	r^{-(d-2)}
	(r^2-x^2)^{(d-3)/2}
	\,
	dr
	\bigg)
	\,
	\mathbb{I}[-1\leq x\leq 1],
	$$
	where~$\mathbb{I}[A]$ denotes the indicator function of~$A$; see~(2.22) in~\cite{Fanetal2000}. Since Part~(i) of the result entails that
	$$
	{\rm BDP}(Q_{\nu_n}(\alpha v))
	=
	\left\lceil  
	n
	\int_\alpha^1 f_1(x)\,dx
	\right\rceil \big/ n
	,
	$$ 
	the result thus follows. 
	(iii) 
	This follows from Part~(i) of the result and from Example~1 in \cite{NagySS}.
	\cqfd  
	\vspace{3mm}

\begin{figure}
		\includegraphics[width=\textwidth]{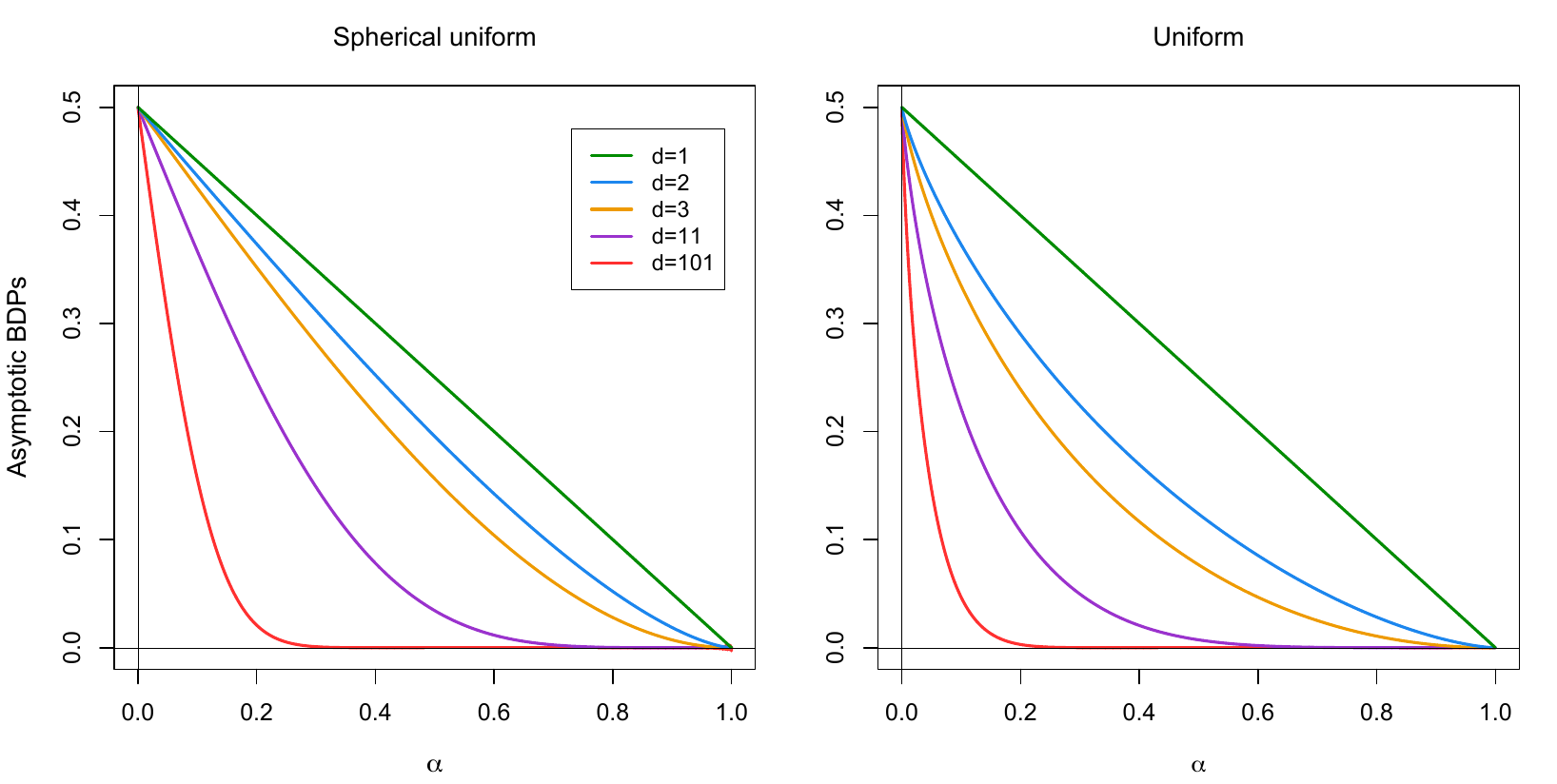}
		\caption{(Left:) Plots, as functions of~$\alpha$, of the asymptotic value of~${\rm BDP}(Q_{\nu_n}(\alpha v))$ obtained for the spherical uniform reference measure in several dimensions~$d$ (this corresponds to the limit of~(\ref{BDPspherunif}) as~$n\to\infty$). (Right:) The corresponding plots for the uniform-on-${\rm cl}(\mathcal{B}_1)$ reference measure (this corresponds to the limit of~(\ref{BDPunif}) as~$n\to\infty$).}
		\label{Fig1}
	\end{figure}

	An illustration of the asymptotic breakdown points associated with~(\ref{BDPspherunif})--(\ref{BDPunif}) is provided in Figure~\ref{Fig1}. For~$d=1$, both asymptotic breakdown points reduce to~$(1-\alpha)/2$, which agrees with the asymptotic breakdown point of univariate quantiles. For both considered reference measures, the asymptotic breakdown point, as expected, is a decreasing function of~$\alpha$ for any fixed dimension~$d$, which is also the case for spatial quantiles (see Corollary~2.2 in \citealp{AIHP}); unlike for spatial quantiles, however, the asymptotic breakdown point of OT quantiles is, still for both considered reference measures, a decreasing function of~$d$ for any fixed order~$\alpha$. Obviously, the asymptotic breakdown point of the OT median (obtained with~$\alpha=0$) is~$1/2$ in any dimension~$d$, which is in line with Corollary~\ref{CorollaryBDPdepth}(ii).

\begin{figure}
		\includegraphics[width=\textwidth]{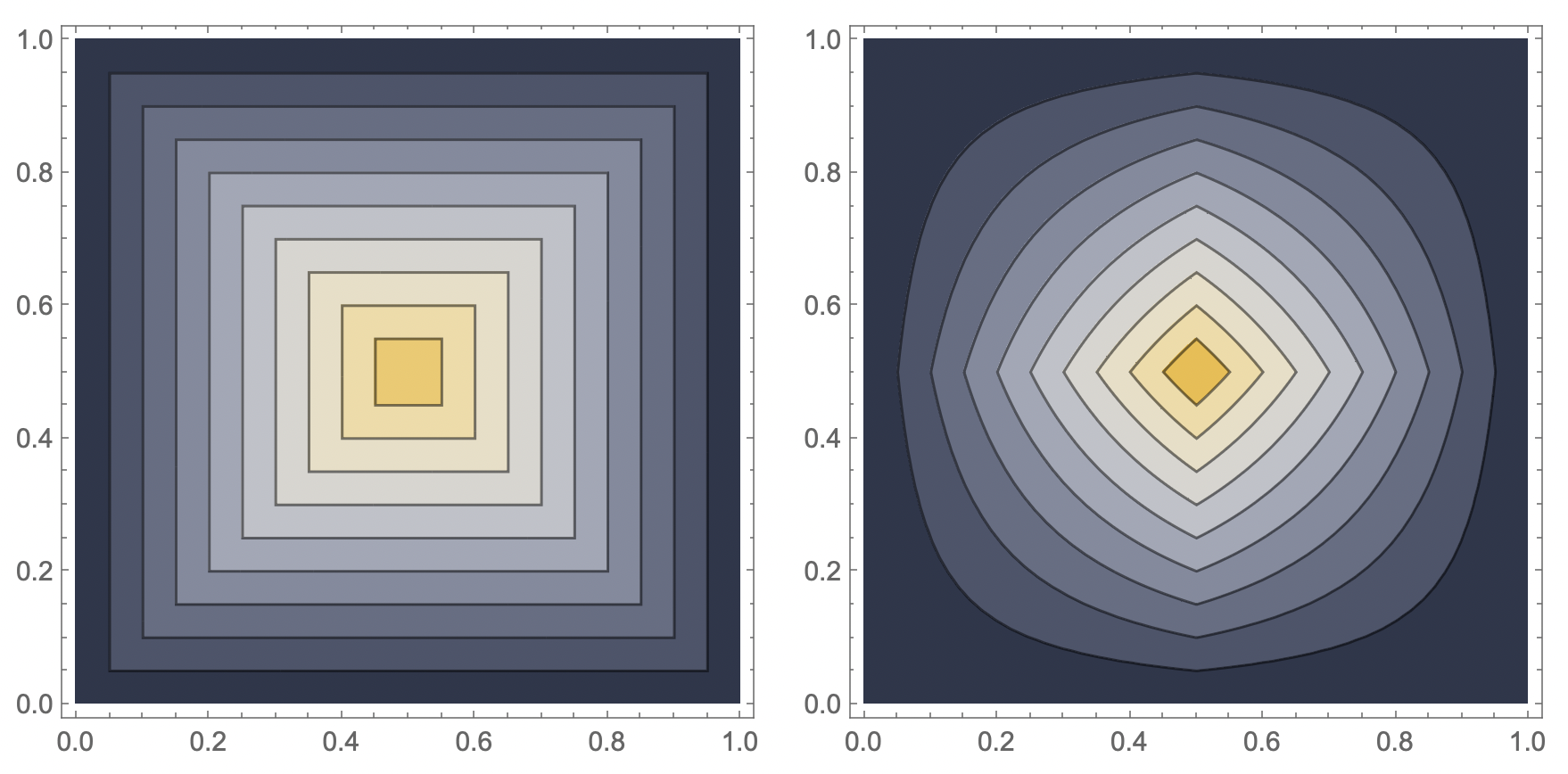}
		\caption{(Left:) The regions used for trimming in the bivariate version of the OT trimmed means in~(\ref{TOtrimmedmeans1}). (Right:) The corresponding (halfspace depth) regions plots used in the bivariate version of the OT trimmed means in~(\ref{TOtrimmedmeans2}).}
		\label{Fig2}
	\end{figure}

	We end this section by commenting on the statistical relevance of our results in the context of robust location estimation. In the univariate case~$d=1$, \emph{$\beta$-trimmed means}, with~$\beta\in[0,1]$, are defined as the average of the~$\lceil n(1-\beta)\rceil$ innermost sample quantiles. The corresponding breakdown point, in the classical Hampel sense, has an asymptotic value equal to~$\beta/2$, which increases from~$0$ to the maximal\footnote{It directly follows from Theorem~2.1 in \cite{LopRou1991} that the asymptotic breakdown point of a translation-equivariant location estimator cannot exceed~$1/2$.} possible value~$1/2$ as one goes from~$\beta=0$ (which provides the sample mean) to~$\beta=1$ (which provides the sample median). In the multivariate setting, a \emph{statistical depth} concept can be used to identify the $\lceil n(1-\beta)\rceil$ most central ``quantiles'', which allows one to consider multivariate~$\beta$-trimmed means. When one adopts a statistical depth function in the sense of \cite{ZuoSer2000A}, this provides affine-equivariant trimmed means, whose breakdown point properties have been studied in the literature; see, e.g.,~\cite{DonGas1992} and \cite{Mas2009} for the halfspace depth, or \cite{Zuo2006} for another celebrated depth, namely the \emph{projection depth}. Now, in order to achieve equivariance under a much broader group of homeomorphic transformations, one may alternatively consider \emph{OT trimmed means}, but this raises the crucial question of how to identify the corresponding innermost quantiles. Focusing on the reference measure~$\mu$ that is put forward in \cite{GS}, namely the uniform measure on~$[0,1]^d$, it is tempting to consider the OT $\beta$-trimmed mean
	\begin{equation}
\label{TOtrimmedmeans1}		
	\frac{
\frac{1}{n}
\sum_{i=1}^n
	x_i \mathbb{I}[ Q_{\nu_n}^{-1}(x_i) \in [\frac{1}{2}-h_\beta,\frac{1}{2}+h_\beta]^d ]
	}
	{
\frac{1}{n}
\sum_{i=1}^n
	\mathbb{I}[ Q_{\nu_n}^{-1}(x_i) \in [\frac{1}{2}-h_\beta,\frac{1}{2}+h_\beta]^d ]
	}
	,
	\end{equation}
	where~$h_\beta$ is the smallest value of~$h(\in[0,\frac{1}{2}])$ such that~$[\frac{1}{2}-h,\frac{1}{2}+h]^d$ contains~$\lceil n(1-\beta)\rceil$ of the OT ``ranks''~$Q_{\nu_n}^{-1}(x_i)$, $i=1,\ldots,n$;  here, $Q_{\nu_n}$ obviously denotes the OT map obtained with the considered reference measure~$\mu$ and the empirical target measure~$\nu_n$ associated with the sample~$x_1,\ldots,x_n$ at hand. In other words, trimming is based on concentric hypercubes in the hypercubical support~$\mathcal{S}=[0,1]^d$ (see the left panel of Figure~\ref{Fig2}). Interestingly, our results reveal that such natural OT $\beta$-trimmed means would actually exhibit a smaller breakdown point than the OT $\beta$-trimmed means defined as
	\begin{equation}
\label{TOtrimmedmeans2}		
	\frac{
\frac{1}{n}
\sum_{i=1}^n
	x_i \mathbb{I}[ Q_{\nu_n}^{-1}(x_i) \in R(h_\beta,\mu) ]
	}
	{
\frac{1}{n}
\sum_{i=1}^n
	\mathbb{I}[ Q_{\nu_n}^{-1}(x_i) \in R(h_\beta,\mu) ]
	}
	,
	\end{equation}
where
$$
R(h,\mu)
:=
\big\{ 
u\in \mathcal{S}: H\!D(u,\mu)\geq h
\big\} 
$$
is the \emph{halfspace depth region} of order~$h$ associated with~$\mu$ and where $h_\beta$ is the largest value of~$h(\in[0,1])$ such that~$R(h,\mu)$ contains~$\lceil n(1-\beta) \rceil$ of the OT ranks~$Q_{\nu_n}^{-1}(x_i)$, $i=1,\ldots,n$ (see the right panel of Figure~\ref{Fig2}). While studying thoroughly the properties of such OT trimmed means is clearly beyond the scope of the present paper, it is remarkable that halfspace depth is relevant even when multivariate trimmed means are based on OT quantiles/ranks rather than on halfspace depth ones.

	\section{Proof of the upper bound}
	\label{secProofUpperBound}

	In this section, we prove the following result, that establishes the upper bound in Theorem~\ref{TheorBDP}(i).

	\begin{prop}
		\label{PropBDPUpperBound}
		Let $\mu$ be an absolutely continuous probability measure with \textcolor{black}{convex} support $\mathcal{S}$ and let~$\nu$ be the discrete probability measure with atoms~$x_1,\ldots,x_n\in\mathbb{R}^d$ and weights~$\lambda_1,\ldots,\lambda_n$. Fix~$u\in \mathcal{S}$. Then, 
		$
		{\rm BDP}(Q_\nu(u))
		\leq
		\min
		\{
		\sum_{i\in I}\lambda_i: \emptyset \neq I\subset\{1,\ldots,n\}
		\ \mathrm{ such\ that }\
		\sum_{i\in I}\lambda_i\geq H\!D(u,\mu) 
		\}
		.
		$
	\end{prop}

	\textcolor{black}{To prove this result, we need the following lemma which will allow us to show that, if~$\sum_{i\in I}\lambda_i\geq H\!D(u,\mu)$, then any neighborhood of~$u$ contains a (fixed) set with positive \mbox{$\mu$-measure} that is contained in the power cell of an arbitrarily strongly contaminated atom} (in the lemma, for~$A\subset\R$, we define~$\sup A$ as~$-\infty$ if~$A=\emptyset$, as~$+\infty$ if~$A$ is not upper-bounded, and as the usual supremum of~$A$ otherwise).
	\vspace{1mm}

	\begin{lem}
		\label{LemContinuityForLowerBound}
		Let~$\mu$ be an absolutely continuous measure over~$\R^d$ with \textcolor{black}{convex} support~$\mathcal{S}$. Fix~$u\in {\rm int}(\mathcal{S})$
		and let~$v_0\in\mathcal{S}^{d-1}$ be such that~$\mathcal{H}_u=\{z\in\R^d: \langle v_0,z-u\rangle \geq 0\}$ is a minimal halfspace in the sense that~$\mu(\mathcal{H}_u)=H\!D(u,\mu)$ $($existence follows from the absolute continuity of~$\mu)$. Then, the map from~$\mathcal{S}^{d-1}$ to $\mathbb{R}\cup\{\pm\infty\}$ defined by
		$$
		v
		\mapsto
		s(v)
		:=
		\sup\big\{ s\in\R : \mu(\mathcal{H}_{u,v,s})\geq H\!D(u,\mu) \big\}
		,
		$$
		where we let
		$
		\mathcal{H}_{u,v,s}
		:=
		\{z\in\R^d: \langle v,z-u \rangle\geq s\}
		$, satisfies the following properties:
\begin{enumerate}
								\item[(i)]	
 $s(\cdot)$ takes its values in~$\R$, 
 								\item[(ii)]	
 $s(\cdot)$ is continuous over~$\mathcal{S}^{d-1}$, 
 								\item[(iii)]	
 $s(v)\geq 0$ for any~$v\in\mathcal{S}^{d-1}$, and 
 								\item[(iv)]	
 $s(v_0)=0$.
\end{enumerate}		
	\end{lem}

	\noindent
	{\sc Proof of Lemma~\ref{LemContinuityForLowerBound}.}
	(i) Fix~$v\in\mathcal{S}^{d-1}$ and let~$A=A_{u,v}:=\{ s\in\R : \mu(\mathcal{H}_{u,v,s})\geq H\!D(u,\mu) \}$. Since~$H\!D(u,\mu)\leq 1/2$ (see, e.g., Lemma~1 from \citealp{RouStr2004}) and $\mu(\mathcal{H}_{u,v,s})\to \mu(\R^d)=1$ as~$s\to-\infty$, the set~$A$ is non-empty. Moreover, since~$H\!D(u,\mu)>0$ (Lemma~\ref{LemHDinSupport}) and $\mu(\mathcal{H}_{u,v,s})\to \mu(\emptyset)=0$ as~$s\to \infty$, there exists~$M\in\R$ such that~$M>s$ for any~$s\in A$. Thus, $s(v)=\sup A\in\R$.
	(ii) Absolute continuity of~$\mu$ implies that the map~$(v,s)\mapsto \mu(\mathcal{H}_{u,v,s})$ is continuous on~$\mathcal{S}^{d-1}\times\R$ \textcolor{black}{(Lemma~\ref{LemContDCT})}. For any~$v$, the map $s\mapsto \mu(\mathcal{H}_{u,v,s})$ is then monotone decreasing and continuous on~$\R$. Thus,
	\begin{equation}
		\label{equivdefsv}
		\mu(\mathcal{H}_{u,v,s(v)})=H\!D(u,\mu)
	\end{equation}
	and there exists~$\eta>0$ such that~$H\!D(u,\mu)/2<\mu(\mathcal{H}_{u,v,s})<3/4$ for any~$s\in (s(v)-\eta,s(v)+\eta)$. \textcolor{black}{The proof proceeds by applying an implicit function theorem to~(\ref{equivdefsv}); since smoothness of~$s\mapsto \mu(\mathcal{H}_{u,v,s})$ is not guaranteed, we will rely on an implicit function theorem for strictly monotone functions (see Theorem~1H.3 from \citealp{Dont2014}).} Would there exist~$s_1,s_2\in(s(v)-\eta,s(v)+\eta)$ with~$s_1<s_2$ and~$\mu(\mathcal{H}_{u,v,s_1})=\mu(\mathcal{H}_{u,v,s_2})$, then
	$$
	\big(\mathcal{S}\cap {\rm cl}(\mathcal{H}_{u,v,s_1}^c) \big)
	\cup 
	\big(\mathcal{S}\cap \mathcal{H}_{u,v,s_2} \big)
	,
	$$  
	would be a \textcolor{black}{strict subset}\footnote{\textcolor{black}{The strict nature of the subset results from the convexity of~$\mathcal{S}$.}} of~$\mathcal{S}$  that is closed and has~$\mu$-measure one, which would contradict the fact that~$\mathcal{S}$ is the support of~$\mu$. Therefore, $s\mapsto \mu(\mathcal{H}_{u,v,s})$ is strictly monotone decreasing in~$(s(v)-\eta,s(v)+\eta)$, so that~$s(v)$ is uniquely defined through~(\ref{equivdefsv}). Since a trivial compactness argument allows one to show that~$\eta>0$ can be chosen independently of~$v$, continuity of~$s(\cdot)$ then follows \textcolor{black}{by applying Theorem~1H.3 from \cite{Dont2014}} to the map~$f(v,s)=-\mu(\mathcal{H}_{u,v,s})+H\!D(u,\mu)$. 
	(iii) From the definition of halfspace depth, $\mu(\mathcal{H}_{u,v,0})\geq H\!D(u,\mu)$ for any~$v\in\mathcal{S}^{d-1}$, which establishes the result.
	(iv) We have seen in the proof of~(ii) that~$s(v_0)$ is uniquely defined through~$\mu(\mathcal{H}_{u,v_0,s(v_0)})=H\!D(u,\mu)$. Since, by assumption, $\mathcal{H}_u=\mathcal{H}_{u,v_0,0}$ is such that~$\mu(\mathcal{H}_{u,v_0,0})=H\!D(u,\mu)$, the result then follows. 
	\cqfd
	\vspace{3mm}


	We can now prove Proposition~\ref{PropBDPUpperBound}.
	\vspace{3mm}

	\noindent
	{\sc Proof of Proposition~\ref{PropBDPUpperBound}.}
	First, we consider the case~$u\in{\rm int}(\mathcal{S})$.
	Fix~$I\subset\{1,\ldots,n\}$ with
	$$
	\sum_{i\in I}\lambda_i\geq H\!D(u,\mu) 
	$$
	(of course, $I$ must then be non-empty). 
	Fix~$v_0\in\mathcal{S}^{d-1}$ such that~$\mathcal{H}_u=\{y\in\R^d: \langle v_0,y-u \rangle \geq 0\}$ is a minimal halfspace (in the same sense as in Lemma~\ref{LemContinuityForLowerBound}) \textcolor{black}{and let~$y_R:=u+Rv_0$ for any~$R>0$. For~$R$ large enough, the contaminated target measure
	\begin{equation}
		\label{tanR}
		\tilde{\nu}_{I,R}
		:=
		\sum_{i\in I^c} \lambda_i \textcolor{black}{\delta_{x_i}}
		+
		\bigg(\sum_{i\in I} \lambda_i \bigg) \textcolor{black}{\delta_{y_R}}
	\end{equation}
	(recall that~$I^c=\{1,\ldots,n\}\setminus I$ collects the indices of the uncontaminated atoms)} has~$m=(\# I^c)+1$ atoms. \textcolor{black}{The proof proceeds by constructing a net~$\{D_\delta\}_{\delta>0}$ of sets with $\mu(D_\delta)>0$ and~$D_\delta\subset\mathcal{B}_\delta(u)$ for all $\delta$ small enough, such that the OT map between $\mu$ and the contaminated target measure~$\tilde{\nu}_{I,R}$, namely $Q_{\tilde{\nu}_{I,R}}$, maps any~$x\in D_\delta$ to~$y_R$. As~$y_R$ will tend to infinity, this will ensure the breakdown; see~(\ref{toccc1})--(\ref{toccc2}) below.}
	
	The power cell of~$y_R$ for the contaminated target measure in~(\ref{tanR})
	is given by
	$$
	{\rm Lag}^{w}_\mathcal{X}(y_R)
	=
	\bigcap_{i\in I^c} (\mathcal{S} \cap \mathcal{H}_{i,R})
	,
	$$ 
	with
	\begin{eqnarray*}
		\mathcal{H}_{i,R}
		& := & 
		\{y\in\mathbb{R}^d:\| y-y_R \|^2-w_{0,R}\leq \| y-x_i\|^2-w_{i,R}\}
		\\[2mm]
		& = & 
		\{y\in\R^d:  \langle v_{i,R} ,y  \rangle \geq s_{i,R}\}
		,
	\end{eqnarray*}
	where $w=(w_{i,R})_{i\in \{0\}\cup I^c}$ is a weight vector\footnote{Here, $w_{0,R}$ is the weight of $y_R$, and $w_{i,R}$, for $i\in I^c$, is the weight of the corresponding atom~$x_i$.} that is adapted to the couple of measures $(\mu, \tilde{\nu}_{I,R})$, where we let $v_{i,R}:=(y_R-x_i)/\|y_R-x_i\|$, and where~$s_{i,R}$ is a real number depending on~$y_R$, $x_i$, $w_{0,R}$ and~$w_{i,R}$.  
	Since the weight vector~$w$ is adapted, we have
	$$
	\mu(\mathcal{H}_{i,R})
	\geq
	\mu({\rm Lag}^{w}_\mathcal{X}(y_R))
	=
	\sum_{i\in I}\lambda_i
	\geq 
	H\!D(u,\mu)
	$$
	for any~$i\in I^c$, so that, using the notation from Lemma~\ref{LemContinuityForLowerBound}, 
	$\mathcal{H}_{i,R}
	\supset
	\mathcal{H}_{u,v_{i,R},s(v_{i,R})}
	$.

	
	Now, fix~$\delta>0$ small enough so that~$\mathcal{B}_\delta(u)\subset\mathcal{S}$. Let~$\eta\in(0,\frac{1}{10})$ be such that 
	$$
	s(v)=|s(v)-s(v_0)|<\frac{\delta}{2} 
	$$ 
	if~$v\in C_{v_0}(\eta):=\{v\in\mathcal{S}^{d-1}:  \langle v_0,v  \rangle \geq 1-\eta\}$ (see Lemma~\ref{LemContinuityForLowerBound})\footnote{\textcolor{black}{This lemma can be used since~$\{ C_v(\eta)\}_{v,\eta}$ is a basis for the topology in~$\mathcal{S}^{d-1}$.}}. For~$R$ large enough, $v_{i,R}\in C_{v_0}(\eta)$ for any~$i\in I^c$. Therefore, for~$R$ large enough,
	$$
	\mathcal{H}_{i,R}
	\supset
	\mathcal{H}_{u,v_{i,R},s(v_{i,R})}
	\supset
	\bigcap_{v\in C_{v_0}(\eta)}
	\ 
	\mathcal{H}_{u,v,s(v)}
	\supset
	\bigcap_{v\in C_{v_0}(\eta)} 
	\ 
	\mathcal{H}_{u,v,\delta/2}
	$$
	for any~$i\in I^c$, so that
	$$
	{\rm Lag}^{w}_\mathcal{X}(y_R)
	=
	\bigcap_{i\in I^c} 
	\
	(\mathcal{S} \cap \mathcal{H}_{i,R})
	\supset
	\bigcap_{v\in C_{v_0}(\eta)}
	(\mathcal{B}_\delta(u) \cap \mathcal{H}_{u,v,\delta/2})
	=:
	D_\delta
	.
	$$
	Note that~$D_\delta$ is a subset of~$\mathcal{B}_\delta(u)$ (hence of~$\mathcal{S}$) that does not depend on~$R$.
	
	\color{black}
	Let us now show that, with~$z_0:=u+3\delta v_0/4$, we have~$\mathcal{B}_{\delta/8}(z_0)\subset D_\delta$. The triangle inequality implies that $\mathcal{B}_{\delta/8}(z_0)\subset \mathcal{B}_\delta(u)$. Now, for any~$v\in C_{v_0}(\eta)$, the distance between~$z_0$ and the boundary hyperplane of~$\mathcal{H}_{u,v,\delta/2}$ is
\begin{equation}
\label{ddDdelta}	
		d(z_0,\partial \mathcal{H}_{u,v,\delta/2})
= \bigg|  \langle v,z_0-u \rangle - \frac{\delta}{2} \bigg|
= \bigg| \frac{3\delta  \langle v,v_0 \rangle }{4} - \frac{\delta}{2} \bigg|
\geq  \frac{3\delta(1-\eta)}{4} - \frac{\delta}{2} 
>  \frac{\delta}{8}
,
\end{equation}
where we used the fact that~$\eta\in (0,\frac{1}{10})$. Since we trivially have~$\mathcal{B}_{\delta/8}(z_0)\subset \mathcal{H}_{u,v_0,\delta/2}$, (\ref{ddDdelta}) implies that~$\mathcal{B}_{\delta/8}(z_0)\subset \mathcal{H}_{u,v,\delta/2}$ for any~$v\in C_{v_0}(\eta)$, so that we indeed have~$\mathcal{B}_{\delta/8}(z_0)\subset D_\delta(\subset\mathcal{S})$. This readily entails that~$\mu(D_\delta)\geq \mu(\mathcal{B}_{\delta/8}(z_0))>0$ (if we would have~$\mu(\mathcal{B}_{\delta/8}(z_0))=0$, then the closed set~$\mathcal{S}\setminus (\mathcal{B}_{\delta/8}(z_0))$ would be a proper subset of~$\mathcal{S}$ with $\mu$-measure one, which would contradict the fact that~$\mathcal{S}$ is the support of~$\mu$). 

With the contaminated target measure~$\tilde{\nu}_{I,R}$ in~(\ref{tanR}), we then have
\color{black}
	\begin{eqnarray}
		\int_{\mathcal{B}_\delta(u)}
		\|Q_\nu(x)-Q_{\tilde{\nu}_{I,R}}(x)\| \, d\mu(x)
		& \geq  & 
		\int_{D_\delta}
		\|Q_\nu(x)-y_R\| \, d\mu(x)
\nonumber
		\\[2mm]
		& \geq  & 
		\big(
		R
        \textcolor{black}{-
        \|u\|}
        -\max(\|x_1\|,\ldots,\|x_n\|)
		\big) 
		\mu(D_\delta)
		.
\label{toccc1}
	\end{eqnarray}
	Since~$\mu(D_\delta)$ is positive and does not depend~$R$, this entails that
	\begin{equation}
	\sup	
	\int_{\mathcal{B}_\delta(u)}
	\|Q_\nu(x)-Q_{\tilde{\nu}_{I}}(x)\| \, d\mu(x)
	\geq
	\lim_{R\to \infty}
	\int_{\mathcal{B}_\delta(u)}
	\|Q_\nu(x)-Q_{\tilde{\nu}_{I,R}}(x)\| \, d\mu(x)
	=
	\infty
	.
\label{toccc2}
	\end{equation}
	Consequently,
	\begin{eqnarray*}
		\lefteqn{
			\hspace{-4mm}
			\Bigg\{
			\sum_{i\in I}\lambda_i: \emptyset \neq I\subset\{1,\ldots,n\}
			\ \, \mathrm{ such\ that }\, \
			\sum_{i\in I}\lambda_i\geq H\!D(u,\mu) 
			\Bigg\}
		}
		\\[2mm]
		& & 
		\hspace{6mm}
		\subset
		\Bigg\{
		\sum_{i\in I}\lambda_i: I\subset\{1,\ldots,n\}
		\ \mathrm{ such\ that } \
		\\[2mm]
		& & 
		\hspace{23mm}
		\sup
		\int_{\mathcal{B}_\delta(u)}
		\|Q_\nu(x)-Q_{\tilde{\nu}_I}(x)\| \, d\mu(x)
		=
		\infty
		\ \
		\forall \delta>0
		\textrm{ small enough}
		\Bigg\}
		.
	\end{eqnarray*}
\textcolor{black}{Since the constraint that~$\delta$ is small enough is actually superfluous, the} result then follows from the definition of the BDP in~(\ref{defBDP}).
	\vspace{3mm}
	
	We now turn to the case~$u\in \mathcal{S}\setminus{\rm int}(\mathcal{S})$. Since~$\mathcal{S}$ is convex and~$u\in \partial \mathcal{S}$, the supporting hyperplane theorem \textcolor{black}{(see Theorem~11.6 in \citealp{Rock1970})} ensures that there exists~$v_0\in\mathcal{S}^{d-1}$ such that~$\mathcal{S}\subset\{y\in\R^d:  \langle v_0,y-u \rangle \leq 0\}$. In particular, $H\!D(u,\mu)=0$ (and~$\mathcal{H}_u=\{y\in\R^d:  \langle v_0,y-u \rangle \geq 0\}$ is a minimal halfspace). Since we then have
	$$
	\min
	\Bigg\{
	\sum_{i\in I}\lambda_i: \emptyset \neq I\subset\{1,\ldots,n\}
	\ \mathrm{ such\ that }\
	\sum_{i\in I}\lambda_i\geq H\!D(u,\mu) 
	\Bigg\}
	=
	\min(\lambda_1,\ldots,\lambda_n)
	,
	$$
	proving Proposition~\ref{PropBDPUpperBound} only requires to show that~${\rm BDP}(Q_\nu(u))\leq \min(\lambda_1,\ldots,\lambda_n)$. By symmetry, it is of course sufficient to show that 
	\begin{equation}
		\label{toconcl31}
		{\rm BDP}(Q_\nu(u))
		\leq 
		\lambda_1
		.
	\end{equation}
	To do so, we let~$I=\{1\}$ (so that $I^c=\{2,\ldots,n\}$) and we adopt a construction that is  similar to the one in the first part of the proof. \textcolor{black}{Letting again~$y_R:=u+Rv_0$ for any~$R>0$, the contaminated target measure}
	\begin{equation}
		\label{tanR2}
		\tilde{\nu}_{I,R}
		=
		\sum_{i\in I^c} \lambda_i \textcolor{black}{\delta_{x_i}}
		+
		\bigg(\sum_{i\in I} \lambda_i\bigg) \textcolor{black}{\delta_{y_R}}
		=
		\lambda_1 \textcolor{black}{\delta_{y_R}}
		+
		\sum_{i=2}^n \lambda_i \textcolor{black}{\delta_{x_i}}
	\end{equation}
	has~$n$ atoms \textcolor{black}{for~$R$ large enough}. 
	\textcolor{black}{As in the case~$u\in {\rm int}(\mathcal{S})$, the proof proceeds by defining a net~$\{D_\delta\}_{\delta>0}$ of sets with $\mu(D_\delta)>0$ and~$D_\delta\subset\mathcal{B}_\delta(u)$ for all $\delta$ small enough, such that~$\tilde{\nu}_{I,R}$ maps any~$x\in D_\delta$ to~$y_R$.}
	Below, ${\rm Lag}^{w}_\mathcal{X}(y_R)$ still denotes the power cell of~$y_R$ for the contaminated target measure in~(\ref{tanR2}), which is again given by
	$$
	{\rm Lag}^{w}_\mathcal{X}(y_R)
	=
	\bigcap_{i\in I^c} (\mathcal{S} \cap \mathcal{H}_{i,R})
	,
	$$ 
	where
	$
	\mathcal{H}_{i,R}
	=
	\{y\in\R^d:  \langle v_{i,R},y \rangle \geq s_{i,R}\}
	$
	involves the same quantities~$v_{i,R}:=(y_R-x_i)/\|y_R-x_i\|$ and~$s_{i,R}$ as in the first part of the proof. 
	
	Now, let~$u_0=u+s_0 v_0$, where~$s_0:=\sup\{s\in \textcolor{black}{\mathbb{R}}: \mu(\mathcal{H}_{u,v_0,s})\geq \lambda_1/2 \}$. Arguing as in the proof of Lemma~\ref{LemContinuityForLowerBound}(ii) and using the same notation as in this lemma, we have
	\begin{equation}
		\label{snzppp}
		\mu(\mathcal{H}_{u_0,v_0,0})
		=
		\mu(\mathcal{H}_{u,v_0,s_0})
		=
		\frac{\lambda_1}{2}
		\cdot
	\end{equation}
	Obviously, $s_0<0$, so that~$u_0\neq u$. Letting~$\delta_0=\|u-u_0\|/2$, 
	$$
	\mathcal{S} \cap \mathcal{B}_{\delta_0}(u) 
	\subset
	\mathcal{S} \cap \mathcal{H}_{u,v_0,s_0}
	.
	$$ 
	Fix then~$\eta>0$ small enough to have both
	$$
	\mu(\mathcal{H}_{u_0,v,0})
	\in
	\Big(
	\frac{\lambda_1}{4}
	,
	\frac{3\lambda_1}{4}
	\Big)
	$$
	for any~$v\in C_{v_0}(\eta)=\{v\in\mathcal{S}^{d-1}: \langle v_0,v \rangle \geq 1-\eta\}$
	and
	$$
	\mathcal{S} \cap \mathcal{B}_{\delta_0}(u) 
	\subset
	\mathcal{S} \cap 
	\bigg(
	\bigcap_{v\in C_{v_0}(\eta)}
	\ 
	\mathcal{H}_{u_0,v,0}
	\bigg)
	$$
	(existence of such an~$\eta$ follows from~(\ref{snzppp}) and from the continuity of the map~$v\mapsto \mu(\mathcal{H}_{u_0,v,0})$). Then, for~$R$ large enough to have $v_{i,R}\in C_{v_0}(\eta)$ for any~$i\in\{2,\ldots,n\}$, we must have
	\begin{equation}
		\label{aaaaz}
		\mathcal{S} \cap
		\mathcal{H}_{i,R}
		\supset
		\mathcal{S} \cap
		\mathcal{H}_{u_0,v_{i,R},0}
		\supset
		\mathcal{S} \cap
		\bigg(
		\bigcap_{v\in C_{v_0}(\eta)}
		\ 
		\textcolor{black}{\mathcal{H}_{u_0,v,0}}
		\bigg)
		\supset
		\mathcal{S} \cap \mathcal{B}_{\delta_0}(u) 
	\end{equation}
	for any~$i\in\{2,\ldots,n\}$ (the first inclusion in~(\ref{aaaaz}) follows from the 
	\color{black}
	facts that $\mathcal{H}_{i,R}$ and $	\mathcal{H}_{u_0,v_{i,R},0}$ share the same outward normal vector~$v_{i,R}$ and 
	\color{black}
	that, since the weight vector~$w$ is adapted, $\mu(\mathcal{H}_{i,R})\geq \mu({\rm Lag}^{w}_\mathcal{X}(y_{R}))=\lambda_1>3\lambda_1/4>\mu(\mathcal{H}_{u_0,v_{i,R},0})$. Therefore, 
	$$
	{\rm Lag}^{w}_\mathcal{X}(y_R)
	=
	\bigcap_{i\in I^c} 
	\
	(\mathcal{S} \cap \mathcal{H}_{i,R})
	\supset
	\mathcal{S} \cap \mathcal{B}_{\delta_0}(u)
	.
	$$
	For any~$\delta\in(0,\delta_0)$, we thus have
	 (with the contaminated target measure~$\tilde{\nu}_{I,R}$ in~(\ref{tanR2}))
	\begin{eqnarray*}
		\int_{\mathcal{B}_\delta(u)}
		\|Q_\nu(x)-Q_{\tilde{\nu}_{I,R}}(x)\| \, d\mu(x)
		& =  & 
		\int_{\mathcal{B}_\delta(u)}
		\|Q_\nu(x)-y_R\| \, d\mu(x)
		\\[2mm]
		& \geq  & 
		\big(
		R
        \textcolor{black}{-\|u\|}
        -\max(\|x_2\|,\ldots,\|x_n\|)
		\big) 
		\mu(\mathcal{B}_\delta(u))
		.
	\end{eqnarray*}
	Since~$\mu(\mathcal{B}_\delta(u))$ is positive and does not depend~$R$, it follows that
	$$
	\sup	
	\int_{\mathcal{B}_\delta(u)}
	\|Q_\nu(x)-Q_{\tilde{\nu}_{I}}(x)\| \, d\mu(x)
	\geq
	\lim_{R\to \infty}
	\int_{\mathcal{B}_\delta(u)}
	\|Q_\nu(x)-Q_{\tilde{\nu}_{I,R}}(x)\| \, d\mu(x)
	=
	\infty
	.
	$$
	Consequently,
	\begin{eqnarray*}
\lefteqn{		
	\lambda_1
	\in
	\Bigg\{
	\sum_{i\in I}\lambda_i: I\subset\{1,\ldots,n\}
	\ \mathrm{ such\ that } \
}
\\[2mm]
& & 
\hspace{16mm}
	\sup
	\int_{\mathcal{B}_\delta(u)}
	\|Q_\nu(x)-Q_{\tilde{\nu}_{I}}(x)\| \, d\mu(x)
	=
	\infty
	\ \
	\forall \delta>0
	\textrm{ small enough}
	\Bigg\}
	.
	\end{eqnarray*}
	As in the first part of the proof, the inequality~(\ref{toconcl31}), hence also the result, then follow from the definition of the BDP in~(\ref{defBDP}).
	\cqfd
	\vspace{3mm}


	\section{Proof of the lower bound}
	\label{secProofLowerBound}

	We turn to the proof of the lower bound in Theorem~\ref{TheorBDP}(i). First note that the definition of the BDP entails that, for any~$u\in\mathcal{S}$,
	$$
	{\rm BDP}(Q_\nu(u))
	\geq 
	\min(\lambda_1,\ldots,\lambda_n)
	=
	\min
	\Bigg\{
	\sum_{i\in I}\lambda_i: \emptyset \neq I\subset\{1,\ldots,n\}
	\ \mathrm{ such\ that }\
	\sum_{i\in I}\lambda_i\geq 
	0
	\Bigg\}
	.
	$$
	For~$u\in \mathcal{S}\setminus {\rm int}(\mathcal{S})$,
	this rewrites
	$$
	{\rm BDP}(Q_\nu(u))
	\geq 
	\min
	\Bigg\{
	\sum_{i\in I}\lambda_i: \emptyset \neq I\subset\{1,\ldots,n\}
	\ \mathrm{ such\ that }\
	\sum_{i\in I}\lambda_i\geq 
	H\!D(u,\mu)
	\Bigg\}
	,
	$$
	so that it is only for~$u\in {\rm int}(\mathcal{S})$ that one needs to prove the lower bound in Theorem~\ref{TheorBDP}(i). In other words, we need to prove the following result.

	\begin{prop}
		\label{PropBDPLowerBound}
		Let $\mu$ be an absolutely continuous probability measure over~$\R^d$ with \textcolor{black}{convex} support~$\mathcal{S}$ and let~$\nu$ be the discrete probability measure with atoms~$x_1,\ldots,x_n\in\mathbb{R}^d$ and weights~$\lambda_1,\ldots,\lambda_n$. Fix~$u\in {\rm int}(\mathcal{S})$. Then, 
		$
		{\rm BDP}(Q_\nu(u))
		\geq
		\min
		\{
		\sum_{i\in I}\lambda_i: 
		\emptyset \neq I\subset\{1,\ldots,n\}
		\ \mathrm{ such\ that }\
		%
		%
		\sum_{i\in I}\lambda_i
		\geq H\!D(u,\mu) 
		\}
		.
		$
	\end{prop}
	
	We will need the following concepts. For any~$j,\ell$, consider the closed halfspace
	\begin{equation}
	    \label{Halfspacemu}
	\mathcal{H}_{j,\ell}
	:=
	\big\{y\in\mathbb{R}^d:\| y-x_j \|^2-w_{j}\leq \| y-x_\ell\|^2-w_\ell\big\}
	. 
	\end{equation}
	For any set of indices~$L$, any~$j\notin L$, \textcolor{black}{and any~$c>0$ large enough so that~$\mathcal{S}_c:=\mathcal{S}\cap {\rm cl}(\mathcal{B}_c)\neq \emptyset$}, we then introduce the set 
	$$
	\textcolor{black}{\mathcal{S}_{c;j,L}}
	:=
	\textcolor{black}{\mathcal{S}_c}
	\cap
	\bigg(
	\cap_{\ell\in L} 
	\mathcal{H}_{j,\ell}
	\bigg)
	.
	$$
	As an intersection of closed convex sets, \textcolor{black}{$\mathcal{S}_{c;j,L}$} is a closed convex set. \textcolor{black}{If~$\partial \mathcal{S}_{c;j,L}\setminus \partial
	\mathcal{S}_c$ is non-empty, then, for any~$z$ in this set,} there exists~$\bar{\ell} \in L$ such that~$z\in \partial \mathcal{H}_{j,\bar{\ell}}$. By definition, we will then say that the atom~$x_{\bar{\ell}}$ is \emph{active} for~$\textcolor{black}{\mathcal{S}_{c;j,L}}$. 
    
    \textcolor{black}{Lemma~\ref{Lem1} below, which will be used to prove Proposition~\ref{PropBDPLowerBound}, analyzes the geometric behavior of power cells. Specifically, given the power cell of an atom located outside $\mathcal{B}_R$ and the union of the power cells of the atoms contained within a smaller ball~$\mathcal{B}_r$, we show that the hyperplanes forming their boundaries---that is, the hyperplanes associated with active atoms---are close to each other. More precisely, within a ball of fixed radius~$c>0$, the maximal distance between such hyperplanes remains uniformly controlled when $R \gg r$. The lemma is essential for understanding the behavior of the power cells of the uncontaminated atoms, which plays a crucial role in the proof of the lower bound, as this proof proceeds by identifying a ball around $u$ that lies entirely within these power cells.
}

\begin{lem}
	\label{Lem1}
	Let~$x_1,\ldots,x_m\in\R^d$ \textcolor{black}{and let $\mu$ be an absolutely continuous probability measure over~$\R^d$ with \textcolor{black}{convex} support~$\mathcal{S}$}. \textcolor{black}{Fix~$c>0$ such that~$\mathcal{S}_c=\mathcal{S}\cap {\rm cl}(\mathcal{B}_c)\neq \emptyset$.}  Fix~$I\subset \{1,\ldots,m\}$ and~$r>0$ such that~$x_i\in \mathcal{B}_r$ for any~$i\in I^c$. Then, for any~$\varepsilon>0$, there exists $R_\varepsilon>r$ (depending only on $\varepsilon$, $c$ and~$r$) such that for any $R\geq R_\varepsilon$, any~$i\in I^c$, any~$j\in I$ for which~$x_j\notin \mathcal{B}_R$, and any~$i_j$ such that~$x_{i_j}$ is active for~$\mathcal{S}_{c;j,I^c}$, we have
	$
	d(y,\partial \mathcal{H}_{i_j,j})
	\leq \varepsilon
	$
	for any $y\in \textcolor{black}{\mathcal{S}_c} \cap \mathcal{H}_{i,j}\cap \mathcal{H}_{j,i_j}$ \textcolor{black}{(all halfspaces involve~$\mu$ through the adapted weight vector~$w$ in~(\ref{Halfspacemu})).} 
\end{lem}

\noindent
{\sc Proof of Lemma~\ref{Lem1}.}
Fix~$\varepsilon>0$, $i\in I^c$, $j\in I$ with~$x_j\notin \mathcal{B}_R$ (at this stage, $R$ is only a fixed number that is strictly larger than~$r$), and~$i_j$ such that~$x_{i_j}$ is active for~$\mathcal{S}_{c;j,I^c}$. If~$x_i=x_{i_j}$, then~$\mathcal{H}_{i,j}\cap \mathcal{H}_{j,i_j}=\partial \mathcal{H}_{i_j,j}$, and we trivially have that~$
d(y,\partial \mathcal{H}_{i_j,j})
=0	\leq \varepsilon
$
for any $y\in \textcolor{black}{\mathcal{S}_c}\cap \mathcal{H}_{i,j}\cap \mathcal{H}_{j,i_j}$. Assume thus that~$x_i\neq x_{i_j}$. 
\emph{Ad absurdum}, take~$y\in \textcolor{black}{\mathcal{S}_c} \cap \mathcal{H}_{i,j}\cap \mathcal{H}_{j,i_j}$ with~$d(y,\partial \mathcal{H}_{i_j,j})> \varepsilon$. Since~$x_{i_j}$ is active for~$\mathcal{S}_{c;j,I^c}$, the set~$\mathcal{S}_{c;j,I^c}\cap \partial\mathcal{H}_{j,i_j}$ is non-empty. Fixing then~$z\in \mathcal{S}_{c;j,I^c}\cap \partial\mathcal{H}_{j,i_j}$ arbitrarily, we have
$$
\| z-x_j\|^2-w_j \leq \| z-x_i\|^2-w_i.
$$
Moreover, since $y\in \mathcal{H}_{i,j}$,  
$$
\| y-x_i\|^2-w_i \leq	\| y-x_j\|^2-w_j
.
$$
Adding up both last expressions provides
\begin{equation*}
	\| z-x_j\|^2-w_j+\| y-x_i\|^2-w_i
	\leq \| z-x_i\|^2-w_i+	\| y-x_j\|^2-w_j
	,
\end{equation*}
or equivalently, 
\begin{equation*}
	\| z-x_j\|^2-\| y-x_j\|^2\leq \| z-x_i\|^2-\| y-x_i\|^2
	.
\end{equation*}
But $\|z-x_i\|^2\leq (c+r)^2$ because $z\in\mathcal{S}_c\subset\mathcal{B}_c$ and $x_i\in\mathcal{B}_r$, and similarly, $\|y-x_i\|^2\leq (c+r)^2$ (since $y\in\textcolor{black}{\mathcal{S}_c}$). Therefore, 
\begin{equation}
	\label{tocontra22}
	\| z-x_j\|^2-\| y-x_j\|^2\leq 2(c+r)^2
	.
\end{equation}

Now, since~$z\in \mathcal{S}_c \cap \partial \mathcal{H}_{j,i_j}$, we have 	
\begin{eqnarray}
	\|z-x_j\|^2-\|y-x_j\|^2
	&   =   &   
	\|z\|^2 - \|y\|^2 + 2  \langle x_j, y-z \rangle 
	\nonumber
	\\[2mm]
	&   =   &   
	\|z\|^2 - \|y\|^2
	+ 2  \langle x_{i_j}, y-z \rangle 
	+ 2  \langle x_j-x_{i_j}, y-z \rangle 
	\nonumber
	\\[2mm]
	&   \geq   &   
	- \|y\|^2
	- 2 \|x_{i_j}\| \|y-z\|
	+ 2  \langle x_j-x_{i_j},y-z \rangle 
	\nonumber
	\\[2mm]
	&   \geq   &   
	- \|y\|^2
	- 2 \|x_{i_j}\| (\|y\|+\|z\|)
	+ 2  \langle x_j-x_{i_j} , y-z \rangle 
	\nonumber
	\\[2mm]
	&   \geq   &   
	- c^2
	- 4cr
	+ 2  \langle x_j-x_{i_j},y-z \rangle 
	,
	\label{pretocontra}
\end{eqnarray}
where we used the fact that~$y,z\in \mathcal{S}_c\subset \mathcal{B}_c$. Since
\begin{eqnarray*}
	\mathcal{H}_{j,i_j}
	& = & 
	\{y\in\mathbb{R}^d:\| y-x_j \|^2-w_j \leq \| y-x_{i_j}\|^2-w_{i_j} \}
	\\[2mm]
	& = & 
	\{y\in\R^d: 2 \langle x_j-x_{i_j},y \rangle  \geq w_{i_j} - w_j + \| x_j\|^2 - \| x_{i_j} \|^2 =: a
	\}
	,
\end{eqnarray*}
using that~$y\in\mathcal{H}_{j,i_j}$ and~$z\in\partial\mathcal{H}_{j,i_j}$ yields
\begin{eqnarray*}
\lefteqn{	
\varepsilon
<
d(y,\partial \mathcal{H}_{i_j,j})
=
\frac{|2 \langle x_j-x_{i_j},y \rangle -a|}{2\|x_j-x_{i_j}\|}
=
\frac{2 \langle x_j-x_{i_j},y \rangle -a}{2\|x_j-x_{i_j}\|}
}
\\[2mm]
& & 
\hspace{41mm}
=
\frac{ \langle x_j-x_{i_j},y-z \rangle }{\|x_j-x_{i_j}\|}
\leq
\frac{ \langle x_j-x_{i_j},y-z \rangle }{R-r}
,
\end{eqnarray*}
so that~(\ref{pretocontra}) provides 
\begin{equation}
	\label{tocontra}	
	\|z-x_j\|^2-\|y-x_j\|^2
	>
	- c^2
	- 4cr
	+ 2\varepsilon (R-r) 
	.
\end{equation}
Therefore, there exists~$R_\varepsilon>r$ such that for any~$R\geq R_\varepsilon$, we have that $\| z-x_j\|^2-\| y-x_j\|^2>2(c+r)^2$, which contradicts~(\ref{tocontra22}). We conclude that, for~$R\geq R_\varepsilon$, any~$y\in \textcolor{black}{\mathcal{S}_c}\cap \mathcal{H}_{i,j}\cap \mathcal{H}_{j,i_j}$ must satisfy~$d(y,\partial \mathcal{H}_{i_j,j})\leq \varepsilon$.
\cqfd
\vspace{3mm}



We can now prove Proposition~\ref{PropBDPLowerBound}.
\vspace{3mm}

\noindent
{\sc Proof of Proposition~\ref{PropBDPLowerBound}.}
\color{black}
We will show that 
\begin{eqnarray}
	\lefteqn{
		\hspace{-4mm}
		\Bigg\{
		\sum_{i\in I}\lambda_i: 		\emptyset \neq I\subset\{1,\ldots,n\}
		\ \, \mathrm{ such\ that }\, \
		\sum_{i\in I}\lambda_i< H\!D(u,\mu) 
		\Bigg\}
	}
	\nonumber
	\\[2mm]
	& & 
	\hspace{6mm}
	\subset
	\Bigg\{
	\sum_{i\in I}\lambda_i: I\subset\{1,\ldots,n\}
	\ \mathrm{ such\ that } \
	\label{toshowlower}
	\\[2mm]
	& & 
	\hspace{23mm}
	\sup
	\int_{\mathcal{B}_\delta(u)}
	\|Q_\nu(x)-Q_{\tilde{\nu}_I}(x)\| \, d\mu(x)
	<
	\infty
	\ \
	\textrm{for some }
    \delta>0
	\Bigg\}
	.
	\nonumber
\end{eqnarray}
Since this implies that
\begin{eqnarray*}
	\lefteqn{
		\hspace{-5mm}
		\Bigg\{
		\sum_{i\in I}\lambda_i: I\subset\{1,\ldots,n\}
		\ \mathrm{ such\ that } \
}
	\\[2mm]
	& & 
	\hspace{22mm}
		\sup
		\int_{\mathcal{B}_\delta(u)}
		\|Q_\nu(x)-Q_{\tilde{\nu}_I}(x)\| \, d\mu(x)
		=
		\infty
		\ \
		\forall \delta>0
		\Bigg\}
	\\[2mm]
	& & 
	\hspace{5mm}
	\subset
	\Bigg\{
	\sum_{i\in I}\lambda_i: 		\emptyset \neq I\subset\{1,\ldots,n\}
	\ \, \mathrm{ such\ that }\, \
	\sum_{i\in I}\lambda_i \geq H\!D(u,\mu) 
	\Bigg\}
	,
	\nonumber
\end{eqnarray*}
the result will follow. 
\color{black}
In order to prove~(\ref{toshowlower}), fix~$I\subset\{1,\ldots,n\}$, $I$ non-empty, with
$$
\sum_{i\in I}\lambda_i
<
H\!D(u,\mu) 
.
$$ 
We need to show that 
there exists~$\delta>0$ such that
$$
\sup
\int_{\mathcal{B}_\delta(u)}
\|Q_\nu(x)-Q_{\tilde{\nu}_I}(x)\| \, d\mu(x)
<
\infty
,
$$
where the supremum is taken over the target measures
$$
\tilde{\nu}_I=\sum_{j=1}^m\tilde{\lambda}_j\delta_{\tilde{x}_j}
,
$$
where the atoms $\tilde{x}_1,\ldots,\tilde{x}_m\in\mathbb{R}^d$ include~$x_i$, $i\in I^c$, and where the weights $\tilde{\lambda}_1,\ldots,\tilde{\lambda}_m$ are so that~$\tilde{\nu}(\{x_i\})=\lambda_i$ for any~$i\in I^c$. It is of course enough to show that there exists~$\delta>0$ such that, for any sequence~$(\tilde{\nu}_\ell)$ of such target measures, the sequence
\begin{equation}
	\label{tib}	
	\Bigg(
	\int_{\mathcal{B}_\delta(u)} 
	\|Q_\nu(x)-Q_{\tilde{\nu}_\ell}(x)\| \, d\mu(x)
	\Bigg)
\end{equation}
is bounded. 

\emph{Ad absurdum}, assume that for any~$\delta>0$ there exists~$(\tilde{\nu}_\ell)$ such that~(\ref{tib})
is unbounded. Fix then~$\delta>0$ small (we will describe later how small) and let~$(\tilde{\nu}_\ell)$ such that~(\ref{tib}) is unbounded. Upon extraction of a subsequence, we may assume that~(\ref{tib}) diverges to infinity. Denoting as~$\tilde{x}_{i,\ell}$, $i=1,\ldots,m_\ell$, the atoms of~$\tilde{\nu}_\ell$, we must then have that the sequence
\begin{equation}
	\label{tib2}
	\bigg( 
	\max_{i=1,\ldots,m_\ell}
	\|\tilde{x}_{i,\ell}\|
	\bigg) 
\end{equation}
is unbounded (otherwise, (\ref{tib}) would be trivially bounded). Upon extraction of a further subsequence, we may assume that~(\ref{tib2}) diverges to infinity. 

\textcolor{black}{Let~$b>0$ such that
$$
\sum_{i\in I}\lambda_i
<
H\!D(u,\mu) 
-
3b
$$
and pick then~$c>0$ large enough so that~$\mathcal{S}_c=\mathcal{S}\cap {\rm cl}(\mathcal{B}_c)$ satisfies~$\mu(\mathcal{S}_c)\geq 1-b$.} 
For any~$\varepsilon\in(0,1)$, let then~$R_\varepsilon$ be as in Lemma~\ref{Lem1} applied with~\textcolor{black}{this~$c$, the subset~$I$ fixed above,} and with~$r>0$ large enough so that the collection of uncontaminated atoms~$\{x_i:i\in I^c\}$ is a subset of~$\mathcal{B}_r$. Denote by~$J\subset \{1,\ldots,m_\ell\}$ the collection of indices~$j$ such that~$\tilde{x}_{j,\ell}\notin \mathcal{B}_{R_\varepsilon}$ (since the sequence in~(\ref{tib2}) diverges to infinity, this collection is non-empty for~$\ell$ large enough). To keep the notation as light as possible, we do not stress dependence on~$\ell$ in~$J$, nor in the quantities/sets we introduce in the rest of the proof.
Consider then the following sets:
\begin{itemize}
	\item The set~$U_\varepsilon$ is defined as the union of the interior of the power cells of the uncontaminated atoms:
	$$
	U_\varepsilon
	:=
	\cup_{i\in I^c} \, {\rm int}({\rm Lag}^{w}_\mathcal{X}(i))
	.
	$$
    \textcolor{black}{Since}
	\begin{equation}
		\label{sfo}
		\mu(U_\varepsilon)
		=
		\mu\Big(
		\!
		\cup_{i\in I^c} 
		{\rm Lag}^{w}_\mathcal{X}(i)
		\Big)
		=
		\sum_{i\in I^c}\lambda_i
		=
		1-\sum_{i\in I}\lambda_i
		>
		1-H\!D(u,\mu)
\textcolor{black}{+3b}
		,
	\end{equation}
    \textcolor{black}{the set
    	$$
	U_{\varepsilon,c}
	:=
{\rm int}(\mathcal{S}_c)\cap	
	U_{\varepsilon}
	$$
    satisfies~$\mu(U_{\varepsilon,c})
\geq
\mu(U_{\varepsilon})-b>1-H\!D(u,\mu)+2b>0$. Thus, 
the open set~$U_{\varepsilon,c}$ is non-empty.} 
	\vspace{2mm}
	\item The set~$Z_\varepsilon$ is defined as
	$$
	Z_\varepsilon
	:= 
	\mathcal{S} \cap 
	\bigg( 
	\cup_{j\in J} 
	\cap_{i\in I^c} 
	\mathcal{H}_{j,i}
	\bigg)
	.
	$$
	Since the power cell of~$\tilde{x}_j(=\tilde{x}_{j,\ell})$ is trivially a subset of $\mathcal{S}\cap (\cap_{i\in I^c} \mathcal{H}_{j,i})$, we have
	$$
	\cup_{j\in J}
    \,
    {\rm Lag}^{w}_\mathcal{X}(j)
	\subset
	\cup_{j\in J} 
	\bigg( 
	\mathcal{S} \cap 
	(\cap_{i\in I^c} 
	\mathcal{H}_{j,i})
	\bigg)
	=
	Z_\varepsilon
	,
	$$
	that is, the union of the power cells of the atoms indexed by~$J$ is a subset of~$Z_\varepsilon$. Similarly, the interior of the power cell of~$\tilde{x}_i$ is a subset of $	\mathcal{S} 
	\cap
	(\cap_{j\in J} 
\,
	 {\rm int}(\mathcal{H}_{i,j}))$, where we used the fact that~${\rm int}(A\cap B)={\rm int}(A)\cap {\rm int}(B)$ for any~$A,B\subset \R^d$, so that
	$
	U_\varepsilon
	\subset
	\mathcal{S} 
	\cap
	( 
	\cup_{i\in I^c} 
	\cap_{j\in J} 
{\rm int}(\mathcal{H}_{i,j})
	)
	$. 
	Therefore,
	\begin{equation}
	    \label{Uepsincll}
U_\varepsilon 
	\subset
	\mathcal{S} 
	\cap
	\bigg( 
	\cap_{j\in J} 
	\cup_{i\in I^c} 
	{\rm int}(\mathcal{H}_{i,j})
	\bigg)
	.
	\end{equation}
    Since
	$$
	Z_\varepsilon
	=
	\mathcal{S} 
	\setminus 
	\bigg( 
	\cup_{j\in J} 
	\cap_{i\in I^c} 
	\mathcal{H}_{j,i}
	\bigg)^c
	=
	\mathcal{S} 
	\setminus 
	\bigg( 
	\cap_{j\in J} 
	\cup_{i\in I^c} 
	{\rm int}(\mathcal{H}_{i,j})
	\bigg)
	,
	$$
	it follows that~$U_\varepsilon \subset \mathcal{S}\setminus Z_\varepsilon$\textcolor{black}{, hence also~$U_{\varepsilon,c} \subset \mathcal{S}_c\setminus Z_\varepsilon$}.  

    \color{black}
Before introducing further sets, we make the following point. Let~$J_0$ be the collection of~$j$'s in~$J$ for which there exists at least one active atom for~$\mathcal{S}_{c;j,I^c}$. We then show that
\begin{equation}
    \label{J0}
\mathcal{S}_c\subset \cup_{i\in I^c}\textnormal{int}(\mathcal{H}_{i,j})
\quad
\textrm{ for all }
j\in J\setminus J_0
.
\end{equation}
To do so, fix~$j\in J\setminus J_0$. Since there is no active atom for~$\mathcal{S}_{c;j,I^c}$, we must have either 
$$
\mathcal{S}_c
\subset
\mathcal{S}\setminus 
\big(
\cap_{i\in I^c} 
\mathcal{H}_{j,i}
\big)
\quad
\Big(
=
\mathcal{S}
\cap
\big(
\cup_{i\in I^c} 
{\rm int}(\mathcal{H}_{i,j})
\big)
\Big)
$$
or 
\begin{equation}
    \label{tskkkkk2}
\mathcal{S}_c
\subset
\mathcal{S}\cap
\Big(
\cap_{i\in I^c} \mathcal{H}_{j,i}
\Big)
,
\end{equation}
and it is therefore sufficient to show that~(\ref{tskkkkk2}) cannot hold. \emph{Ad absurdum}, assume that~(\ref{tskkkkk2}) holds. Then, (\ref{Uepsincll}) yields
$$
U_\varepsilon
\subset
\mathcal{S}
\cap
\Big(
\cup_{i\in I^c} 
{\rm int}(\mathcal{H}_{i,j})
\Big)
=
\mathcal{S}
\cap
\Big(
\cap_{i\in I^c} 
\mathcal{H}_{j,i}
\Big)^c
,
$$ 
which provides 
$$
\mu(
\mathcal{S}\cap 
(\cap_{i\in I^c} \mathcal{H}_{j,i})
)
\leq 
\mu(\mathcal{S}\setminus U_\varepsilon)
.
$$
Recalling that~$\mu(U_\varepsilon)>1-H\!D(u,\mu)+3b$ and $\mu(S_c)\geq 1-b$, (\ref{tskkkkk2}) entails that
		$$
        1-b\leq\mu(S_c)
        \leq \mu(\mathcal{S}\cap
        (\cap_{i\in I^c} \mathcal{H}_{j,i})
        )
        \leq  \mu(\mathcal{S}\setminus U_\varepsilon)< H\!D(u,\mu)-3b
        .
        $$
	Since this contradicts the fact that~$H\!D(u,\mu)\leq 1$, the claim~(\ref{J0}) is proved.
\color{black}

	\vspace{2mm}
	\item The sets~$V_\varepsilon$ and~$W_\varepsilon$ are defined as
	\begin{equation}
		\label{defVW}
		V_\varepsilon
		:=
		\textcolor{black}{\mathcal{S}_c} \cap 
		\bigg( 
		\cap_{j\in \textcolor{black}{J_0}} 
		\mathcal{H}_{i_j,j}
		\bigg) 
		\quad
		\textrm{and}
		\quad
		W_\varepsilon
		:=
		\textcolor{black}{\mathcal{S}_c} \cap 
		\bigg( 
		\cap_{j\in \textcolor{black}{J_0}} 
		\mathcal{H}_{i_j,j}^{\varepsilon}
		\bigg) 
		,
	\end{equation}
	respectively, where, for any~$j\in \textcolor{black}{J_0}$, we took~$\tilde{x}_{i_j}$ arbitrarily in the \textcolor{black}{(non-empty)} collection of atoms that are active for~$\textcolor{black}{\mathcal{S}_{c;j,I^c}}$, and where we let~$A^\varepsilon:=\{x\in \R^d: d(x,A)\leq \varepsilon\}$\textcolor{black}{; in~(\ref{defVW}), an intersection over an empty collection of subsets of~$\R^d$ is defined as being~$\R^d$, so that when~$J_0$ is empty, we simply have~$V_\varepsilon=W_\varepsilon=\mathcal{S}_c$. Since~(\ref{J0}) yields that}
	$$
	{\rm int}(V_\varepsilon)
	=
	{\rm int}(\textcolor{black}{\mathcal{S}_c})
	\cap 
	\bigg( 
	\cap_{j\in \textcolor{black}{J_0}} 
	{\rm int}(\mathcal{H}_{i_j,j})
	\bigg) 
	\subset
	\mathcal{S} \cap 
	\bigg( 
	\cap_{j\in J} 
	\cup_{i\in I^c} 
	{\rm int}(\mathcal{H}_{i,j})
	\bigg) 
	,
	$$
	we have that~${\rm int}(V_\varepsilon)\subset S\setminus Z_\varepsilon$. 
    
    \textcolor{black}{Now, assume that $J_0$ is not empty.} Lemma~\ref{Lem1} entails that, for any~$j\in \textcolor{black}{J_0}$, 
	\begin{equation}
		\label{inclusionA}
	\textcolor{black}{\mathcal{S}_c}
	\cap 
	\bigg(
	\cup_{i\in I^c} 
	\mathcal{H}_{i,j}
	\bigg)
\textcolor{black}{\,\cap\ \mathcal{H}_{j,i_j}}
	\subset
	\textcolor{black}{\mathcal{S}_c} 
	\cap 
	\mathcal{H}_{i_j,j}^{\varepsilon}
	.
	\end{equation}
\color{black}
	Since~$\mathcal{S}_c 
	\cap \mathcal{H}_{j,i_j}^c
	=
	\mathcal{S}_c 
	\cap {\rm int}(\mathcal{H}_{i_j,j})
	\subset
	\mathcal{S}_c 
	\cap \mathcal{H}_{i_j,j}^{\varepsilon}$ for any~$j\in \textcolor{black}{J_0}$,
    we also have that, for any~$j\in \textcolor{black}{J_0}$, 
	\begin{equation}
		\label{inclusionB}
	\textcolor{black}{\mathcal{S}_c} 
	\cap 
	\bigg(
	\cup_{i\in I^c} 
	\mathcal{H}_{i,j}
	\bigg)
	\cap \mathcal{H}_{j,i_j}^c
	\subset
	\textcolor{black}{\mathcal{S}_c} 
	\cap 
	\mathcal{H}_{i_j,j}^{\varepsilon}
.
	\end{equation}
Clearly, (\ref{inclusionA})--(\ref{inclusionB}) imply that, for any~$j\in \textcolor{black}{J_0}$,
\color{black}
\begin{equation*}
    	\textcolor{black}{\mathcal{S}_c} 
	\cap 
	\bigg(
	\cup_{i\in I^c} 
	\mathcal{H}_{i,j}
	\bigg)
	\subset
	\textcolor{black}{\mathcal{S}_c} 
	\cap 
	\mathcal{H}_{i_j,j}^{\varepsilon}
.
\end{equation*}
\textcolor{black}{Thus, when $J_0$ is not empty, this  yields}
	$$
\textcolor{black}{U_{\varepsilon,c}}
	\subset
\textcolor{black}{\mathcal{S}_c}
	\setminus
	Z_\varepsilon
	=
\textcolor{black}{\mathcal{S}_c}	\cap 
	\bigg( 
	\cap_{j\in J} 
	\cup_{i\in I^c} 
	{\rm int}(\mathcal{H}_{i,j})
	\bigg)
	\subset
	\textcolor{black}{\mathcal{S}_c}
    \cap 
	\bigg(
	\cap_{j\in \textcolor{black}{J_0}} 
	\cup_{i\in I^c}
	\mathcal{H}_{i,j}
	\bigg)
	\subset
	W_\varepsilon
	.
	$$
\textcolor{black}{Note that the inclusion $U_{\varepsilon,c}\subset W_\varepsilon$ holds as well when~$J_0$ is empty. Indeed, when~$J_0$ is empty, (\ref{J0}) yields
	$$
\textcolor{black}{U_{\varepsilon,c}}
	\subset
\textcolor{black}{\mathcal{S}_c}
	\setminus
	Z_\varepsilon
	=
\textcolor{black}{\mathcal{S}_c}	\cap 
	\bigg( 
	\cap_{j\in J} 
	\cup_{i\in I^c} 
	{\rm int}(\mathcal{H}_{i,j})
	\bigg)
	=
	\textcolor{black}{\mathcal{S}_c}
    =
	W_\varepsilon
	.
	$$}
    
	%
	%
\end{itemize}

Since $\textcolor{black}{U_{\varepsilon,c}}\subset W_\varepsilon$, it follows from~(\ref{sfo}) that 
\color{black} 
\begin{equation}
	\label{UU}
	\mu(W_\varepsilon)
	\geq
	\mu(\textcolor{black}{U_{\varepsilon,c}})
	>
	1-H\!D(u,\mu)
	+2b
	>
	0
\end{equation}
for any~$\varepsilon\in(0,1)$
\color{black}
(which obviously guarantees that~$W_\varepsilon$ is non-empty). Lemma~\ref{LemAbsCont}(i) entails that there exists~$v>0$ (not depending on~$\varepsilon$) such that~$\mathcal{L}(W_\varepsilon)\geq v$ for any~$\varepsilon\in(0,1)$. Now, \color{black} since the convex set~$W_\varepsilon$ is  included~$\mathcal{B}_c$, hence has diameter\footnote{Here, the diameter of a subset~$A$ of~$\R^d$ is defined as~$\sup_{x,y\in A} \|x-y\|$.} at most~$2c$.  \color{black} Therefore, Lemma~\ref{LemmaPolytopeBall} implies that there exists $\zeta\in(0,1)$ such that $W_\varepsilon$ contains a ball of radius~$\zeta$ (note that this radius only depends on~$v$ and~$c$, hence in particular does not depend on~$\varepsilon$ nor on the shape of~$W_\varepsilon$). For any~$\varepsilon\in(0,\zeta/4)$, the definitions of~$V_\varepsilon$ and~$W_\varepsilon$ in~(\ref{defVW}) then entail that $V_\varepsilon$ contains a ball of radius~$\zeta/2$, which guarantees in particular that~$V_\varepsilon$ is non-empty.

Denoting as~$r_{V_\varepsilon}(\geq \zeta/2)$ the radius of the largest closed ball contained in $V_\varepsilon$, Point~(ii) in Lemma A6.5 in \cite{Kal} then yields 
\begin{eqnarray*}
	\mathcal{L}(W_\varepsilon\setminus V_\varepsilon)
	&   \leq   &  
	2
	\big\{
	(1+\varepsilon/r_{V_\varepsilon})^d-1 
	\big\}
	\mathcal{L}(V_\varepsilon)
	\\[2mm]
	&   \leq   &  
	2
	\big\{
	(1+2\varepsilon/\zeta)^d-1 
	\big\}
	\textcolor{black}{\mathcal{L}(\mathcal{S}_c)}
	,
\end{eqnarray*}
where we used the fact that~$V_\varepsilon$ is a subset of~$\textcolor{black}{\mathcal{S}_c}$. Therefore, Lemma~\ref{LemAbsCont}(ii) implies that~$\mu(W_\varepsilon\setminus V_\varepsilon)\to 0$ as~$\varepsilon\to 0$. \color{black} Jointly with~(\ref{UU}), this implies that there exists~$\varepsilon'>0$ such that
\begin{eqnarray*}
	\mu(V_\varepsilon)
	&   =   &  
	\mu(W_\varepsilon)
	-
	\mu(W_\varepsilon\setminus V_\varepsilon)
	\\[2mm]
	&   \geq    &  
	1-H\!D(u,\mu)
	+
	b
\end{eqnarray*}
for any~$\varepsilon\in(0,\varepsilon')$
\color{black}
(again, this ensures that~$V_\varepsilon$ has a non-empty interior).
%
%
%
%
%
%
%
%
%
%
We now use this to prove that there exists~$\delta>0$ such that $\mathcal{B}_\delta(u) \subset {\rm int}(V_\varepsilon)$ for any~$\varepsilon\in(0,\varepsilon')$. To do so, fix $\varepsilon\in(0,\varepsilon')$ arbitrarily. Note first that if~$u\notin {\rm int}(V_\varepsilon)$, then~${\rm int}(V_\varepsilon)$ is a convex subset of~$\mathcal{S}$ that does not contain~$u$ and satisfies~$\mu({\rm int}(V_\varepsilon))=\mu(V_\varepsilon)>1-H\!D(u,\mu)$, which contradicts Lemma~\ref{LemmaLargeConvexSubset}. Thus, $u\in {\rm int}(V_\varepsilon)$. \color{black}
Fix then $\delta>0$ such that $\mu(C)<b/2$ for any band~$C \in \mathcal{C}_{\delta}$, where the notation~$\mathcal{C}_{\delta}$ is introduced in Lemma~\ref{Lemmainfband} (existence of~$\delta$ follows from this lemma). Assume \emph{ad absurdum} that we do not have~$\mathcal{B}_\delta(u) \subset {\rm int}(V_\varepsilon)$. Then, there exists~$z_\varepsilon\in \partial V_\varepsilon$ with
\begin{equation}
\label{diszepsVeps}    
\|u-z_\varepsilon\|
<
\delta
.
\end{equation}
Let~$H$ be an arbitrary hyperplane supporting $V_\varepsilon$ at $z_\varepsilon$ and denote by $H_u$ the hyperplane obtained by translating~$H$ to~$u$. Due to~(\ref{diszepsVeps}), the width of the closed band~$C$ between~$H_u$ and~$H$ is smaller than~$\delta$, hence satisfies $\mu(C)<b/2$. Therefore,
$$
1-H\!D(u,\mu)
+
b
\leq
\mu(V_\varepsilon)
=
\mu(V_\varepsilon\cap C)
+
\mu(V_\varepsilon\cap C^c)
<
\frac{b}{2}
+
\mu(V_\varepsilon\cap C^c)
,
$$
so that~$\mu(V_\varepsilon\cap C^c)>1-H\!D(u,\mu)+(b/2)$. It follows that $V_\varepsilon\cap C^c$ is a convex set (as the intersection between two convex sets) that does not contain~$u$ and whose~$\mu$-measure is strictly larger than~$1-H\!D(u,\mu)$. Since this contradicts Lemma~\ref{LemmaLargeConvexSubset} again, we must have that~$\mathcal{B}_\delta(u) \subset {\rm int}(V_\varepsilon)$. 
\color{black}
Since $\varepsilon\in(0,\varepsilon')$ was fixed arbitrarily and~$\delta$ does not depend on~$\varepsilon$, we proved that~$\mathcal{B}_\delta(u) \subset {\rm int}(V_\varepsilon)$ for any~$\varepsilon\in(0,\varepsilon')$. 

We can now conclude the proof. Fix indeed~$\delta>0$ such that~$\mathcal{B}_\delta(u) \subset {\rm int}(V_\varepsilon)$ for any~$\varepsilon\in(0,\varepsilon')$. As explained at the beginning of the proof, the contradiction argument provides the existence of a sequence of contaminated target measures~$(\tilde{\nu}_\ell)$ such that
\begin{equation}
	\label{tnzmpp}
	\Bigg(
	\int_{\mathcal{B}_\delta(u)} 
	\|Q_\nu(x)-Q_{\tilde{\nu}_\ell}(x)\| \, d\mu(x)
	\Bigg)
	\to 
	\infty
	.
\end{equation}
As explained again at the beginning of the proof, the maximal norm of the atoms of~$(\tilde{\nu}_\ell)$ diverges to infinity as~$\ell$ does, so that, for~$\ell$ large enough, there is at least one atom of~$(\tilde{\nu}_\ell)$ outside~$\mathcal{B}_{R_{\varepsilon'/2}}$. However, since any~$x$ in~$\mathcal{B}_\delta(u)$ then belongs to the power cells of an atom in~$\mathcal{B}_{R_{\varepsilon'/2}}$ (recall indeed that, for any~$\varepsilon\in(0,1)$, we have that~${\rm int}(V_\varepsilon)\subset \mathcal{S}\setminus Z_\varepsilon$ and that the union of the power cell of the atoms outside~$\mathcal{B}_{R_\varepsilon}$ is a subset of~$Z_\varepsilon$), a direct consequence is that
$$
\int_{\mathcal{B}_\delta(u)}
\|Q_\nu(x)-Q_{\tilde{\nu}_\ell}(x)\| \, d\mu(x)
\leq
\int_{\mathcal{B}_\delta(u)}
( \|Q_\nu(x)\| + R_{\varepsilon'/2}) \, d\mu(x)
<
\infty
$$
for any~$\ell$ large enough. Since this contradicts~(\ref{tnzmpp}), the result is proved.
\cqfd


\section{\textcolor{black}{Final comments} and perspectives for future research}
\label{secPersp}

Our work provides the first quantitative analysis of the robustness of (semi-discrete) OT. In particular, it precisely characterizes the finite-sample breakdown point of the OT quantiles from \cite{GS}, and it does so for a very broad class of reference measures~$\mu$. Our results reveal that this breakdown point depends on the geometry of~$\mu$ through the well-known concept of halfspace depth. In particular, it is only for (angularly) symmetric reference measures that the asymptotic breakdown point of the OT median coincides with the one of the univariate median, that is, is equal to~$1/2$. In the context of robust location estimation, our results provide a subtle insight on how to perform multivariate trimming when constructing OT trimmed means.

\textcolor{black}{As just mentioned, our breakdown point results focus on the semi-discrete OT quantiles from \cite{GS}. Remarkably, an interesting independent work, \cite{Ave2024}, has obtained the corresponding results in two complementary frameworks, namely (i) when both the reference and target measures are discrete (which provides the OT quantiles from \citealp{Cheetal2017} and \citealp{Hal21}) and (ii) when both are continuous. While the Tukey halfspace depth is relevant in all frameworks, it is interesting to note that the dimension plays a stronger role in the discrete case than in the semi-discrete and continuous ones. In particular, the discrete OT median from in \cite{Cheetal2017} and \cite{Hal21} seems to have a very low breakdown point in high-dimensions,\footnote{\textcolor{black}{We refer to Remark~3.3 in \cite{Ave2024}, that states that when the support of the discrete reference measure at hand is included in a hyperplane of~$\R^d$ (which is obviously always the case in the high-dimensional framework where~$d\geq n$), then the breakdown point of the OT map is~$1/n$.}} whereas their semi-discrete and continuous analogs will keep breakdown point~$1/2$ at any high-dimensional, angularly symmetric, reference measure. 
As we explain below, this suggests that the semi-discrete approach considered in this paper may provide good results in infinite-dimensional spaces, too.}

Perspectives for future research are rich and diverse. \textcolor{black}{We mention three of them here:}
\begin{itemize}
    \item[(i)]
    Inspection of the proof of the upper bound in Proposition~\ref{PropBDPUpperBound} shows that \textcolor{black}{breakdown of~$Q_\nu(u)$ occurs when putting the whole contamination budget arbitrarily far on a halfline with direction~$v_0$ originating from~$u$, where~$v_0$ is the inner-normal unit vector of an arbitrary (Tukey) minimal halfspace at~$u$} (see the contaminated target measures in~(\ref{tanR}) and~(\ref{tanR2})). This worst-case scenario is thus similar to the one that provides the breakdown point of spatial quantiles; see Theorems~3.2--3.3 in \cite{AIHP}. It would be interesting to see how robust are OT quantiles under other types of contamination, for instance when contamination occurs in another direction~$v$ or when contamination is placed symmetrically in two opposite directions. For spatial quantiles, the resulting breakdown points are strictly higher than the global breakdown point associated with the worst-case scenario above; see Section~4 in \cite{AIHP}. It may be challenging, however, to obtain results of this type for OT quantiles. 
\vspace{2mm}
    \item[(ii)]
\textcolor{black}{As mentioned in the introduction, the optimal transport maps considered in this work are defined only $\mu$-almost everywhere. A direct corollary is that one cannot consider pointwise breakdown points. Addressing this issue by restricting to~$u$'s for which both~$Q_\nu(u)$ and $Q_{\tilde{\nu}_I}(u)$ are well-defined would have considerably reduced the generality of our results. This is why we rather adopted the local averaging around~$u$ in~(\ref{defBDP}) when defining the breakdown point. As pointed out by a Reviewer, this issue already materializes in dimension~$d=1$, where, however, there is a common agreement to consider right-continuous quantile functions. In higher dimensions, \cite{Segers2022} introduced set-valued quantile functions in order to overcome the lack of a canonical choice in~$\mathbb{R}^d$. Since these  are defined everywhere in the support of the reference measure~$\mu$, it would be natural to try and determine the breakdown point of such set-valued semi-discrete quantile functions.}
\vspace{2mm}
    \item[(iii)] 
 The results of this paper focus on finite-dimensional Euclidean spaces. But more and more research efforts are dedicated to defining OT quantiles in infinite-dimensional Hilbert spaces; too; see, e.g., \cite{GHS23}. It would thus be interesting to extend our robustness results to this more general framework. \textcolor{black}{Obviously, the fact that the dimension does not play a role in our results---in contrast to the discrete results in \cite{Ave2024}---is encouraging in this respect. Actually, one of the Reviewers reported that the arguments used in our proofs could indeed be adapted to the infinite-dimensional case: in this framework,  the role of the supporting hyperplane theorem will be played by the Hahn--Banach theorem, whereas suitable assumptions on the measure~$\mu$ will ensure that it does not give mass to the boundary of the cells (this will be the case, e.g., for non-degenerated Gaussian measures).}
\end{itemize}

\textcolor{black}{The important technical challenges these questions raise explain that these are left for future research work.}


\begin{appendix}
\section*{Auxiliary results}
\label{secAppendix} 



%

In this appendix, we prove the results that were used in the proofs of Sections~\ref{secProofUpperBound}--\ref{secProofLowerBound}.
\vspace{1mm}

\begin{lem}
	\label{LemHDinSupport}
	\textcolor{black}{Let~$\mu$ be a measure}  over~$\R^d$ with support~$\mathcal{S}$. Then, $H\!D(u,\mu)>0$ for any~$u\in {\rm int}(\mathcal{S})$. 	
\end{lem}

\noindent
{\sc Proof of Lemma~\ref{LemHDinSupport}.}
\color{black}
\emph{Ad absurdum}, let~$u\in {\rm int}(\mathcal{S})$ such that~$H\!D(u,\mu)=0$. Thus, there exists a sequence~$(H_n)$ of halfspaces such that~$u\in H_n$ for all~$n$ and such that~$\mu(H_n)\to 0$. Obviously, there is no loss of generality to assume that~$u\in \partial H_n$ for all~$n$ (for any~$n$ for which we have~$u\notin \partial H_n$, one may replace~$H_n$ with the closed halfspace that is contained in~$H_n$ and has~$u$ on its boundary hyperplane). Each~$H_n$ is thus of the form~$
		\mathcal{H}_{u,v_n}$ for some~$v_n\in\mathcal{S}^{d-1}$, where we wrote~$\mathcal{H}_{u,v}
		:=
		\{z\in\R^d: \langle v,z-u \rangle\geq 0\}$. Compactness of~$\mathcal{S}^{d-1}$ implies that there exists a subsequence~$(v_{n_\ell})$ that converges in~$\mathcal{S}^{d-1}$, to~$v_0$ say. Fatou's lemma then entails that
		\begin{eqnarray*}
\lefteqn{			
		\mu({\rm int}(H_{u,v_0}))
		\leq
		\int_{\R^d} 
		\liminf_{\ell\to\infty}
		\mathbb{I}[z\in H_{n_\ell},z\in {\rm int}(H_{u,v_0})]\,d\mu(z) 
}
\\[2mm]
& & 
\hspace{3mm}
		\leq
		\int_{\R^d} 
		\liminf_{\ell\to\infty}
		\mathbb{I}[z\in H_{n_\ell}]\,d\mu(z) 
		\leq
		\liminf_{\ell\to\infty}
		\int_{\R^d} \mathbb{I}[z\in H_{n_\ell}]\,d\mu(z) 
		=
		\liminf_{\ell\to\infty}
		\mu(H_{n_\ell}) 
		=
		0
		.
		\end{eqnarray*}
Thus, $\mu({\rm int}(H_{u,v_0}))=0$. Since~$u\in{\rm int}(\mathcal{S})$, there exists~$\delta>0$ such that~$\mathcal{B}_\delta(u)\subset \mathcal{S}$. Thus, $\mathcal{B}_{\delta/2}(u+(\delta/2)v_0)\subset \mathcal{S}\cap {\rm int}(H_{u,v_0})$. We must have that~$\mu(\mathcal{B}_{\delta/2}(u+(\delta/2)v_0))>0$ (otherwise, $\mathcal{S}\setminus \mathcal{B}_{\delta/2}(u+(\delta/2)v_0)$ is a proper closed subset of~$\mathcal{S}$ with $\mu$-measure, which contradicts that~$\mathcal{S}$ is the support of~$\mu$). Therefore,
$$
		\mu({\rm int}(H_{u,v_0}))
\geq
\mu(\mathcal{B}_{\delta/2}(u+(\delta/2)v_0))
>
0
,
$$ 
a contradiction. The result is proved.
%
\color{black}
\cqfd
\vspace{3mm}


\color{black}
\begin{lem}
	\label{LemContDCT}
	Let~$\mu$ be an absolutely continuous measure over~$\R^d$. Fix~$u\in\R^d$, $v\in\mathcal{S}^{d-1}$ and~$s\in\R$. Let~$(v_\ell,s_\ell)$ be a sequence in~$\mathcal{S}^{d-1}\times \R$ that converges to~$(v,s)$.   Then, using the notation introduced in Lemma~\ref{LemContinuityForLowerBound}, $\mu(\mathcal{H}_{u,v_\ell,s_\ell})$ converges to~$\mu(\mathcal{H}_{u,v,s})$. 	
\end{lem}

\noindent
{\sc Proof of Lemma~\ref{LemContDCT}.}
Since absolute continuity of~$\mu$ entails that~$\mu(\partial \mathcal{H}_{u,v,s})=0$, we have
\begin{eqnarray*}
| \mu(\mathcal{H}_{u,v_\ell,s_\ell}) - \mu(\mathcal{H}_{u,v,s}) |
&\leq &
\int_{\R^d}
|
\mathbb{I}_{\mathcal{H}_{u,v_\ell,s_\ell}}(x) 
-
\mathbb{I}_{\mathcal{H}_{u,v,s}}(x) 
|
\,
d\mu(x)
\\[2mm]
&= & 
\int_{\R^d\setminus \partial \mathcal{H}_{u,v,s}}
|
\mathbb{I}_{\mathcal{H}_{u,v_\ell,s_\ell}}(x) 
-
\mathbb{I}_{\mathcal{H}_{u,v,s}}(x) 
|
\,
d\mu(x)
.
\end{eqnarray*}
Since~$(v_\ell,s_\ell)$ converges to~$(v,s)$, we have that $\mathbb{I}_{\mathcal{H}_{u,v_\ell,s_\ell}}(x)$ converges to~$\mathbb{I}_{\mathcal{H}_{u,v,s}}(x)$ for any~$x\in \R^d\setminus \partial \mathcal{H}_{u,v,s}$. A routine application of Lebesgue's Dominated Convergence Theorem thus yields that~$\mu(\mathcal{H}_{u,v_\ell,s_\ell})$ converges to~$\mu(\mathcal{H}_{u,v,s})$. 
\cqfd
\vspace{3mm}

\color{black}


\begin{lem}
	\label{LemAbsCont}
	Let~$\mu$ be a measure over~$\mathbb{R}^d$ that is absolutely continuous with respect to the Lebesgue measure~$\mathcal{L}$. Then, (i) for any~$a>0$,
	$
	\inf
	\{ \mathcal{L}(A): A\in \mathcal{A}_{\mu,a} \}
	>
	0
	$,
	where~$\mathcal{A}_{\mu,a}$ denotes the collection of Borel sets of~$\R^d$ such that $\mu(A)\geq a$. 
	(ii) For any sequence~$(A_n)$ of Borel sets in~$\mathbb{R}^d$ such that~$\mathcal{L}(A_n)\to 0$, we have~$\mu(A_n)\to 0$.
\end{lem}

\noindent
{\sc Proof of Lemma~\ref{LemAbsCont}.}
(i) Fix~$a>0$. Denoting as~$g$ the Lebesgue density of~$\mu$ and as~$\mathbb{I}[A]$ the indicator function of~$A$, the Dominated Convergence Theorem entails that
$$
\int_{\mathbb{R}^d}
g(x) \mathbb{I}[g(x)> k]
\,
d\mathcal{L}(x)
\to 	
0
$$
as~$k$ diverges to infinity, so that there exists~$k_a>0$ for which	
$$
\int_{\mathbb{R}^d}
g(x) \mathbb{I}[g(x)>k_a]
\,
d\mathcal{L}(x)
\leq
\frac{a}{2}
\cdot
$$
For any~$A\in\mathcal{A}_{\mu,a}$, we then have
\begin{equation*}
	a
	\leq
	\mu(A)
	=
	\int_{A} g(x) \mathbb{I}[g(x)> k_a]
	\,
	d\mathcal{L}(x)
	+
	\int_{A} g(x) \mathbb{I}[g(x)\leq k_a]
	\,
	d\mathcal{L}(x)
	\leq
	\frac{a}{2}
	+
	k_a	\mathcal{L}(A)
	,
\end{equation*}
hence also~$\mathcal{L}(A)\geq a/(2k_a)$.
(ii)
Fix~$\varepsilon>0$ and a sequence~$(A_n)$ of Borel sets in~$\mathbb{R}^d$ such that~$\mathcal{L}(A_n)\to 0$. With the same notation as above, we have
\begin{equation*}
	\mu(A_n)
	=
	\int_{A_n} g(x) \mathbb{I}[g(x)> k_\varepsilon]
	\,
	d\mathcal{L}(x)
	+
	\int_{A_n} g(x) \mathbb{I}[g(x)\leq k_\varepsilon]
	\,
	d\mathcal{L}(x)
	\leq
	\frac{\varepsilon}{2}
	+
	k_\varepsilon	\mathcal{L}(A_n)
	,
\end{equation*}
for any~$n$. Therefore, there exists a natural number~$N$ such that~$\mu(A_n)<\varepsilon$ for any~$n\geq N$, which establishes the result.
\cqfd

%
%
%


\begin{lem}
	\label{LemmaPolytopeBall}
	For $r,v>0$, let~$\mathcal{D}_{r,v}$ be the collection of convex sets in~$\R^d$ whose diameter is smaller than or equal to~$r$ and whose volume is larger than or equal to~$v$. Fix~$r$ and $v$ such that~$\mathcal{D}_{r,v}$ is not void. Then, there exists~$\zeta=\zeta_{r,v}>0$ such that any~$D\in \mathcal{D}_{r,v}$ contains a ball of radius~$\zeta$. 
\end{lem} 

\noindent
{\sc Proof of Lemma~\ref{LemmaPolytopeBall}.}
\emph{Ad absurdum}, assume that there exists a sequence of sets~$(D_n)$ in~$\mathcal{D}_{r,v}$ such that~$D_n$ does not contain a ball of radius~$1/n$. Clearly,~$D_n$ is contained in a closed hyper-cylinder of radius~$r/2$ and length~$2/n$. It follows that, for any~$n$ large enough, $D_n$ has a volume that is strictly smaller than~$v$, a contradiction.
\cqfd


\begin{lem}
	\label{LemmaLargeConvexSubset}
	Let~$\mathcal{S}$ be a convex set in~$\R^d$ and let~$\mu$ be an absolutely continuous measure whose support is~$\mathcal{S}$. Fix~$u\in\mathcal{S}$. Then, the maximal $\mu$-measure that can be achieved by a convex subset of~$\mathcal{S}$ that does not contain~$u$ is~$1-H\!D(u,\mu)$. 
\end{lem}

\noindent
{\sc Proof of Lemma~\ref{LemmaLargeConvexSubset}.}
Let~$\mathcal{H}_u^-$ be an arbitrary closed halfspace that contains~$u$ on its boundary hyperplane and that satisfies 
\begin{equation}
	\label{HDpre}
	H\!D(u,\mu)
	=
	\mu(\mathcal{S}\cap \mathcal{H}_u^-)
\end{equation}
(existence follows from absolute continuity of~$\mu$). Thus, $D_u:=\mathcal{S}\cap (\mathcal{H}_u^-)^c$ is a convex subset of~$\mathcal{S}$ that does not contain~$u$, and its $\mu$-measure is~$1-H\!D(u,\mu)$ (this follows from~(\ref{HDpre}) since~$\mu(\mathcal{S}\cap \mathcal{H}_u^-)+\mu(\mathcal{S}\cap (\mathcal{H}_u^-)^c)=\mu(\mathcal{S})=1$). We need to prove that this is indeed the maximal $\mu$-measure that can be achieved by a convex subset of~$\mathcal{S}$ that does not contain~$u$. \emph{Ad absurdum}, assume then that there exists a convex subset,~$D$ say, of~$\mathcal{S}$ that does not contain~$u$ and that satisfies $\mu(D)>1-H\!D(u,\mu)$. Since~$D$ is convex and~$u\notin D$, there exists a hyperplane, $\tilde{\mathcal{H}}_u$ say, containing~$u$ that does not intersect~$D$. Among both closed halfspaces determined by~$\tilde{\mathcal{H}}_u$, denote then as~$\tilde{\mathcal{H}}_u^-$ the one that does not contain~$D$. Since
$
\mu((\tilde{\mathcal{H}}_u^-)^c)
\geq 
\mu(D)
>
1-H\!D(u,\mu)
$,
the closed halfspace~$\tilde{\mathcal{H}}_u^-$ has~$u$ on its boundary hyperplane and satisfies
$
\mu(\tilde{\mathcal{H}}_u^-)
<
H\!D(u,\mu)
$, 
a contradiction to the definition of $H\!D(u,\mu)$ (see the statement of Theorem~\ref{TheorBDP}). 
\cqfd




\begin{lem}
	\label{Lemmainfband}
    \color{black}
	Let~$\mu$ be a measure over~$\mathbb{R}^d$ that is bounded and is absolutely continuous with respect to the Lebesgue measure. For any~$v\in\mathcal{S}^{d-1}$, $c\geq 0$, and~$h>0$, define the band
	$$
	C_{v,c,h}
	:=
	\{x\in\R^d: c\leq  \langle v,x \rangle \leq c+h\}
	$$
	and denote as~$\mathcal{C}_{\delta}:=\{	C_{v,c,h} : v\in\mathcal{S}^{d-1},c\geq 0 , h \in (0,\delta] \}$ the collection of such bands with width at most~$\delta$. Then, for any~$a>0$, there exists~$\delta>0$ such that
	$
	\mu(C)<a
	$
	for any~$C\in \mathcal{C}_{\delta}$. 
\color{black}
\end{lem} 

\noindent
{\sc Proof of Lemma~\ref{Lemmainfband}.}
\color{black}
Without any loss of generality, we may assume that~$a<\mu(\R^d)$. 
Let~$R>0$ be large enough to have
$
\mu(\R^d\setminus\mathcal{B}_R)
<
a
$. Since~$\|x\| \geq |\langle v,x\rangle|\geq c$ for any~$x\in C_{v,c,h}$, we then have that
\begin{equation}
    \label{sjfpp}
\mu(C_{v,c,h})
<
a
\qquad
\forall
v\in\mathcal{S}^{d-1},
\
c\geq R,
\
h>0
.
\end{equation}
The same routine application of the Lebesgue DCT as in the proof of Lemma~\ref{LemContDCT} establishes that
$
(v,c,h)
\mapsto
\mu(C_{v,c,h})
$
is continuous over~$\mathcal{S}^{d-1}\times[0,R]\times [0,1]$. Thus, the sequence~$(g_\ell)$ of functions defined on~$\mathcal{S}^{d-1}\times[0,R]$ by 
$$
g_\ell
:
\mathcal{S}^{d-1}\times[0,R]
\to
\R
:
(v,c) 
\mapsto 
\mu(C_{v,c,1/\ell})
$$
converges pointwise to the zero function (absolute continuity indeed entails that~$\mu(C_{v,c,0})$ for any~$v,c$). Compactness of~$\mathcal{S}^{d-1}\times[0,R]$ implies that~$(g_\ell)$ converges also uniformly to the zero function, so that there exists a positive integer~$N$ such that
$$
\sup
\big\{
\mu(C_{v,c,1/n})
:
v\in\mathcal{S}^{d-1},
\
c\leq R
\big\}
<
a
$$
for any~$n\geq N$. Recalling~(\ref{sjfpp}), we conclude that~$\mu(C)<a$	for any~$C\in \mathcal{C}_{1/N}$.  
\color{black}
\cqfd

\end{appendix}

%
%
%

\begin{acks}[Acknowledgments]
		Davy Paindaveine is also affiliated at the Toulouse School of Economics, Université Toulouse 1 Capitole.
\textcolor{black}{The authors would like to thank the anonymous referees, an Associate Editor and the Editors for their insightful comments, that improved significantly the quality of the paper.}
\end{acks}

\begin{funding}
	Davy Paindaveine was supported by the ``Projet de Recherche'' T.0230.24 from the FNRS (Fonds National pour la Recherche Scientifique), Communauté Française de Belgique, and by a grant from the Fonds Thelam, King Baudouin Foundation. 
\end{funding}

\bibliographystyle{imsart-nameyear} 
\bibliography{Paper}       

\providecommand{\noopsort}[1]{}
\begin{thebibliography}{83}

\bibitem[\protect\citeauthoryear{\'{A}lvarez Esteban
  et~al.}{2008}]{Alvarez2008}
\begin{barticle}[author]
\bauthor{\bparticle{\'{A}lvarez} \bsnm{Esteban},~\bfnm{P.~C.}\binits{P.~C.}},
  \bauthor{\bsnm{Del~Barrio},~\bfnm{E.}\binits{E.}},
  \bauthor{\bsnm{Cuesta-Albertos},~\bfnm{J.}\binits{J.}} \AND
  \bauthor{\bsnm{Matr\'{a}n},~\bfnm{C.}\binits{C.}}
(\byear{2008}).
\btitle{Trimmed comparison of distributions}.
\bjournal{J. Amer. Statist. Assoc.}
\bvolume{103}
\bpages{697--704}.
\end{barticle}
\endbibitem

\bibitem[\protect\citeauthoryear{Ambrosio, Semola and
  Bru\'{e}}{2021}]{Ambrosio2021}
\begin{bbook}[author]
\bauthor{\bsnm{Ambrosio},~\bfnm{L.}\binits{L.}},
  \bauthor{\bsnm{Semola},~\bfnm{D.}\binits{D.}} \AND
  \bauthor{\bsnm{Bru\'{e}},~\bfnm{E.}\binits{E.}}
(\byear{2021}).
\btitle{Lectures on Optimal Transport}.
\bpublisher{Springer}, \baddress{Cham, Switzerland}.
\end{bbook}
\endbibitem

\bibitem[\protect\citeauthoryear{Aurenhammer}{1987}]{AHA87}
\begin{barticle}[author]
\bauthor{\bsnm{Aurenhammer},~\bfnm{F.}\binits{F.}}
(\byear{1987}).
\btitle{Power diagrams: properties, algorithms and applications}.
\bjournal{SIAM J. Comput.}
\bvolume{16}
\bpages{78--96}.
\end{barticle}
\endbibitem

\bibitem[\protect\citeauthoryear{Aurenhammer, Hoffmann and
  Aronov}{1998}]{AHA98}
\begin{barticle}[author]
\bauthor{\bsnm{Aurenhammer},~\bfnm{F.}\binits{F.}},
  \bauthor{\bsnm{Hoffmann},~\bfnm{F.}\binits{F.}} \AND
  \bauthor{\bsnm{Aronov},~\bfnm{B.}\binits{B.}}
(\byear{1998}).
\btitle{Minkowski-type theorems and least-squares clustering}.
\bjournal{Algorithmica}
\bvolume{20}
\bpages{61--76}.
\end{barticle}
\endbibitem

\bibitem[\protect\citeauthoryear{Avella-Medina and
  Gonz\'{a}lez-Sanz}{2024}]{Ave2024}
\begin{barticle}[author]
\bauthor{\bsnm{Avella-Medina},~\bfnm{M.}\binits{M.}} \AND
  \bauthor{\bsnm{Gonz\'{a}lez-Sanz},~\bfnm{A.}\binits{A.}}
(\byear{2024}).
\btitle{On the breakdown point of transport-based quantiles}.
\bjournal{arXiv:2410.1655}.
\end{barticle}
\endbibitem

\bibitem[\protect\citeauthoryear{Balaji, Chellappa and
  Feizi}{2020}]{Balayi2020}
\begin{barticle}[author]
\bauthor{\bsnm{Balaji},~\bfnm{Y.}\binits{Y.}},
  \bauthor{\bsnm{Chellappa},~\bfnm{R.}\binits{R.}} \AND
  \bauthor{\bsnm{Feizi},~\bfnm{S.}\binits{S.}}
(\byear{2020}).
\btitle{Robust optimal transport with applications in generative modeling and
  domain adaptation}.
\bjournal{Adv. Neural Inf. Process. Syst.}
\bvolume{33}
\bpages{12934--12944}.
\end{barticle}
\endbibitem

\bibitem[\protect\citeauthoryear{Bansil and Kitagawa}{2022}]{Bansil}
\begin{barticle}[author]
\bauthor{\bsnm{Bansil},~\bfnm{M.}\binits{M.}} \AND
  \bauthor{\bsnm{Kitagawa},~\bfnm{J.}\binits{J.}}
(\byear{2022}).
\btitle{Quantitative stability in the geometry of semi-discrete optimal
  transport}.
\bjournal{Int. Math. Res. Not. IMRN}
\bvolume{2022}
\bpages{7354--7389}.
\end{barticle}
\endbibitem

\bibitem[\protect\citeauthoryear{Bassetti, Bodini and
  Regazzini}{2006}]{Bassetti2006}
\begin{barticle}[author]
\bauthor{\bsnm{Bassetti},~\bfnm{F.}\binits{F.}},
  \bauthor{\bsnm{Bodini},~\bfnm{A.}\binits{A.}} \AND
  \bauthor{\bsnm{Regazzini},~\bfnm{E.}\binits{E.}}
(\byear{2006}).
\btitle{On minimum Kantorovich distance estimators}.
\bjournal{Statist. Probab. Lett.}
\bvolume{76}
\bpages{1298--1302}.
\end{barticle}
\endbibitem

\bibitem[\protect\citeauthoryear{Bercu and Bigot}{2021}]{Bercu2021}
\begin{barticle}[author]
\bauthor{\bsnm{Bercu},~\bfnm{B.}\binits{B.}} \AND
  \bauthor{\bsnm{Bigot},~\bfnm{J.}\binits{J.}}
(\byear{2021}).
\btitle{Asymptotic distribution and convergence rates of stochastic algorithms
  for entropic optimal transportation between probability measures}.
\bjournal{Ann. Statist.}
\bvolume{49}
\bpages{968--987}.
\end{barticle}
\endbibitem

\bibitem[\protect\citeauthoryear{Bercu, Bigot and Thurin}{2024}]{BigotDirect24}
\begin{barticle}[author]
\bauthor{\bsnm{Bercu},~\bfnm{B.}\binits{B.}},
  \bauthor{\bsnm{Bigot},~\bfnm{J.}\binits{J.}} \AND
  \bauthor{\bsnm{Thurin},~\bfnm{G.}\binits{G.}}
(\byear{2024}).
\btitle{Regularized estimation of Monge-Kantorovich quantiles for spherical
  data}.
\bjournal{arXiv:2407.02085}.
\end{barticle}
\endbibitem

\bibitem[\protect\citeauthoryear{Boissard and Le~Gouic}{2014}]{Boissard2014}
\begin{barticle}[author]
\bauthor{\bsnm{Boissard},~\bfnm{E.}\binits{E.}} \AND
  \bauthor{\bsnm{Le~Gouic},~\bfnm{T.}\binits{T.}}
(\byear{2014}).
\btitle{On the mean speed of convergence of empirical and occupation measures
  in Wasserstein distance}.
\bjournal{Ann. Inst. Henri Poincar\'{e} Probab. Stat.}
\bvolume{50}
\bpages{539--563}.
\end{barticle}
\endbibitem

\bibitem[\protect\citeauthoryear{Brenier}{1991}]{Brenier}
\begin{barticle}[author]
\bauthor{\bsnm{Brenier},~\bfnm{Y.}\binits{Y.}}
(\byear{1991}).
\btitle{Polar factorization and monotone rearrangement of vector-valued
  functions}.
\bjournal{Comm. Pure Appl. Math.}
\bvolume{44}
\bpages{375--417}.
\end{barticle}
\endbibitem

\bibitem[\protect\citeauthoryear{Chaudhuri}{1996}]{Cha1996}
\begin{barticle}[author]
\bauthor{\bsnm{Chaudhuri},~\bfnm{Probal}\binits{P.}}
(\byear{1996}).
\btitle{On a geometric notion of quantiles for multivariate data}.
\bjournal{J. Amer. Statist. Assoc.}
\bvolume{91}
\bpages{862--872}.
\end{barticle}
\endbibitem

\bibitem[\protect\citeauthoryear{Chernozhukov et~al.}{2017}]{Cheetal2017}
\begin{barticle}[author]
\bauthor{\bsnm{Chernozhukov},~\bfnm{Victor}\binits{V.}},
  \bauthor{\bsnm{Galichon},~\bfnm{Alfred}\binits{A.}},
  \bauthor{\bsnm{Hallin},~\bfnm{Marc}\binits{M.}} \AND
  \bauthor{\bsnm{Henry},~\bfnm{Marc}\binits{M.}}
(\byear{2017}).
\btitle{Monge--Kantorovich depth, quantiles, ranks and signs}.
\bjournal{Ann. Statist.}
\bvolume{45}
\bpages{223--256}.
\end{barticle}
\endbibitem

\bibitem[\protect\citeauthoryear{Deb, Bhattacharya and Sen}{2021}]{Bha2021}
\begin{barticle}[author]
\bauthor{\bsnm{Deb},~\bfnm{N.}\binits{N.}},
  \bauthor{\bsnm{Bhattacharya},~\bfnm{B.~B.}\binits{B.~B.}} \AND
  \bauthor{\bsnm{Sen},~\bfnm{B.}\binits{B.}}
(\byear{2021}).
\btitle{Pitman efficiency lower bounds for multivariate distribution-free tests
  based on optimal transport}.
\bjournal{arXiv:2104.01986}.
\end{barticle}
\endbibitem

\bibitem[\protect\citeauthoryear{Deb and Sen}{2023}]{DebSen2022}
\begin{barticle}[author]
\bauthor{\bsnm{Deb},~\bfnm{N.}\binits{N.}} \AND
  \bauthor{\bsnm{Sen},~\bfnm{B.}\binits{B.}}
(\byear{2023}).
\btitle{Multivariate rank-based distribution-free nonparametric testing using
  measure transportation}.
\bjournal{J. Amer. Statist. Assoc.}
\bvolume{118}
\bpages{192--207}.
\end{barticle}
\endbibitem

\bibitem[\protect\citeauthoryear{del Barrio, Gonz\'{a}lez~Sanz and
  Loubes}{2024}]{delBarrio2024}
\begin{barticle}[author]
\bauthor{\bparticle{del} \bsnm{Barrio},~\bfnm{E.}\binits{E.}},
  \bauthor{\bsnm{Gonz\'{a}lez~Sanz},~\bfnm{A.}\binits{A.}} \AND
  \bauthor{\bsnm{Loubes},~\bfnm{J.~M.}\binits{J.~M.}}
(\byear{2024}).
\btitle{Central limit theorems for semi-discrete Wasserstein distances}.
\bjournal{Bernoulli}
\bvolume{30}
\bpages{554--580}.
\end{barticle}
\endbibitem

\bibitem[\protect\citeauthoryear{del Barrio and Loubes}{2019}]{delBarrio2019}
\begin{barticle}[author]
\bauthor{\bparticle{del} \bsnm{Barrio},~\bfnm{E.}\binits{E.}} \AND
  \bauthor{\bsnm{Loubes},~\bfnm{J.~M.}\binits{J.~M.}}
(\byear{2019}).
\btitle{Central limit theorems for empirical transportation cost in general
  dimension}.
\bjournal{Ann. Probab.}
\bvolume{47}
\bpages{926--951}.
\end{barticle}
\endbibitem

\bibitem[\protect\citeauthoryear{Donoho and Gasko}{1992}]{DonGas1992}
\begin{barticle}[author]
\bauthor{\bsnm{Donoho},~\bfnm{David~L.}\binits{D.~L.}} \AND
  \bauthor{\bsnm{Gasko},~\bfnm{Miriam}\binits{M.}}
(\byear{1992}).
\btitle{Breakdown properties of location estimates based on halfspace depth and
  projected outlyingness}.
\bjournal{Ann. Statist.}
\bvolume{20}
\bpages{1803--1827}.
\end{barticle}
\endbibitem

\bibitem[\protect\citeauthoryear{Donoho and Huber}{1983}]{Dono83}
\begin{bincollection}[author]
\bauthor{\bsnm{Donoho},~\bfnm{David~L.}\binits{D.~L.}} \AND
  \bauthor{\bsnm{Huber},~\bfnm{Peter~J.}\binits{P.~J.}}
(\byear{1983}).
\btitle{The notion of breakdown point}.
In \bbooktitle{A Festschrift for Erich Lehmann}
(\beditor{\bfnm{P.~J.}\binits{P.~J.}~\bsnm{Bickel}},
  \beditor{\bfnm{K.}\binits{K.}~\bsnm{Doksum}} \AND
  \beditor{\bfnm{J.~L.}\binits{J.~L.}~\bsnm{Hodges}}, eds.)
\bpages{157--184}.
\bpublisher{Wadsworth}, \baddress{Belmont, CA}.
\end{bincollection}
\endbibitem

\bibitem[\protect\citeauthoryear{Dontchev and Rockafellar}{2014}]{Dont2014}
\begin{bbook}[author]
\bauthor{\bsnm{Dontchev},~\bfnm{A.~L.}\binits{A.~L.}} \AND
  \bauthor{\bsnm{Rockafellar},~\bfnm{R.~T.}\binits{R.~T.}}
(\byear{2014}).
\btitle{Implicit Functions and Solution Mappings. A View From Variational
  Analysis},
\bedition{2nd} ed.
\bpublisher{Springer}, \baddress{New York}.
\end{bbook}
\endbibitem

\bibitem[\protect\citeauthoryear{Dudley and Koltchinski}{1992}]{DuKol92}
\begin{bunpublished}[author]
\bauthor{\bsnm{Dudley},~\bfnm{R.~M.}\binits{R.~M.}} \AND
  \bauthor{\bsnm{Koltchinski},~\bfnm{V.~I.}\binits{V.~I.}}
(\byear{1992}).
\btitle{The spatial quantiles}.
\bnote{Unpublished manuscript.}
\end{bunpublished}
\endbibitem

\bibitem[\protect\citeauthoryear{Fang, Kotz and Ng}{1990}]{Fanetal2000}
\begin{bbook}[author]
\bauthor{\bsnm{Fang},~\bfnm{Kai-Tai}\binits{K.-T.}},
  \bauthor{\bsnm{Kotz},~\bfnm{Samuel}\binits{S.}} \AND
  \bauthor{\bsnm{Ng},~\bfnm{Kai-Wang}\binits{K.-W.}}
(\byear{1990}).
\btitle{Symmetric Multivariate and Related Distributions}.
\bpublisher{Springer Science \& Business Media}, \baddress{New Delhi}.
\end{bbook}
\endbibitem

\bibitem[\protect\citeauthoryear{Galichon}{2016}]{Galichon2016}
\begin{bbook}[author]
\bauthor{\bsnm{Galichon},~\bfnm{A.}\binits{A.}}
(\byear{2016}).
\btitle{Optimal Transport Methods in Economics}.
\bpublisher{Princeton University Press}.
\end{bbook}
\endbibitem

\bibitem[\protect\citeauthoryear{Ghosal and Sen}{2022a}]{GS}
\begin{barticle}[author]
\bauthor{\bsnm{Ghosal},~\bfnm{P.}\binits{P.}} \AND
  \bauthor{\bsnm{Sen},~\bfnm{B.}\binits{B.}}
(\byear{2022}a).
\btitle{Multivariate ranks and quantiles using optimal transport: Consistency,
  rates and nonparametric testing}.
\bjournal{Ann. Statist.}
\bvolume{50}
\bpages{1012--1037}.
\end{barticle}
\endbibitem

\bibitem[\protect\citeauthoryear{Ghosal and Sen}{2022b}]{GSsupp}
\begin{barticle}[author]
\bauthor{\bsnm{Ghosal},~\bfnm{P.}\binits{P.}} \AND
  \bauthor{\bsnm{Sen},~\bfnm{B.}\binits{B.}}
(\byear{2022}b).
\btitle{Supplement to ``Multivariate ranks and quantiles using optimal
  transport: Consistency, rates and nonparametric testing''}.
\bjournal{Ann. Statist.}
\bvolume{50}
\bpages{1012--1037}.
\end{barticle}
\endbibitem

\bibitem[\protect\citeauthoryear{Goldfeld et~al.}{2024}]{Goldfeld2024}
\begin{barticle}[author]
\bauthor{\bsnm{Goldfeld},~\bfnm{Z.}\binits{Z.}},
  \bauthor{\bsnm{Kato},~\bfnm{K.}\binits{K.}},
  \bauthor{\bsnm{Rioux},~\bfnm{G.}\binits{G.}} \AND
  \bauthor{\bsnm{Sadhu},~\bfnm{R.}\binits{R.}}
(\byear{2024}).
\btitle{Limit theorems for entropic optimal transport maps and Sinkhorn
  divergence}.
\bjournal{Electron. J. Stat.}
\bvolume{18}
\bpages{980--1041}.
\end{barticle}
\endbibitem

\bibitem[\protect\citeauthoryear{Gonz\'{a}lez-Sanz, Hallin and
  Sen}{2025}]{GHS23}
\begin{barticle}[author]
\bauthor{\bsnm{Gonz\'{a}lez-Sanz},~\bfnm{A.}\binits{A.}},
  \bauthor{\bsnm{Hallin},~\bfnm{M.}\binits{M.}} \AND
  \bauthor{\bsnm{Sen},~\bfnm{B.}\binits{B.}}
(\byear{2025}).
\btitle{Monotone measure-transportation maps in Hilbert spaces,
  \textcolor{black}{with statistical applications}}.
\bjournal{\textcolor{black}{Bernoulli},}
\bpages{\textcolor{black}{to appear}}.
\end{barticle}
\endbibitem

\bibitem[\protect\citeauthoryear{Gonz\'{a}lez~Sanz, Loubes and
  Niles-Weed}{2022}]{Gonz2022}
\begin{barticle}[author]
\bauthor{\bsnm{Gonz\'{a}lez~Sanz},~\bfnm{A.}\binits{A.}},
  \bauthor{\bsnm{Loubes},~\bfnm{J.~M.}\binits{J.~M.}} \AND
  \bauthor{\bsnm{Niles-Weed},~\bfnm{Jonathan}\binits{J.}}
(\byear{2022}).
\btitle{Weak limits of entropy regularized optimal transport; potentials, plans
  and divergences}.
\bjournal{arXiv:2207.07427}.
\end{barticle}
\endbibitem

\bibitem[\protect\citeauthoryear{Gonz\'{a}lez-Sanz and
  Sheng}{2024}]{GonSheng2024}
\begin{barticle}[author]
\bauthor{\bsnm{Gonz\'{a}lez-Sanz},~\bfnm{A.}\binits{A.}} \AND
  \bauthor{\bsnm{Sheng},~\bfnm{S.}\binits{S.}}
(\byear{2024}).
\btitle{Linearization of Monge--Amp\`{e}re equations and statistical
  applications}.
\bjournal{arXiv:2408.06534}.
\end{barticle}
\endbibitem

\bibitem[\protect\citeauthoryear{Hallin, D. and Liu}{2022}]{HallinJASA22}
\begin{barticle}[author]
\bauthor{\bsnm{Hallin},~\bfnm{M.}\binits{M.}},
  \bauthor{\bsnm{D.},~\bfnm{La~Vecchia}\binits{L.~V.}} \AND
  \bauthor{\bsnm{Liu},~\bfnm{H.}\binits{H.}}
(\byear{2022}).
\btitle{Center-outward R-estimation for semi-parametric VARMA models}.
\bjournal{J. Amer. Statist. Assoc.}
\bvolume{117}
\bpages{925--938}.
\end{barticle}
\endbibitem

\bibitem[\protect\citeauthoryear{Hallin, Hlubinka and
  Hudecov\'{a}}{2023}]{HallinHlu2023}
\begin{barticle}[author]
\bauthor{\bsnm{Hallin},~\bfnm{M.}\binits{M.}},
  \bauthor{\bsnm{Hlubinka},~\bfnm{Daniel}\binits{D.}} \AND
  \bauthor{\bsnm{Hudecov\'{a}},~\bfnm{{\v{S}}}\binits{{\v{S}}.}}
(\byear{2023}).
\btitle{Efficient fully distribution-free center-outward rank tests for
  multiple-output regression and MANOVA}.
\bjournal{J. Amer. Statist. Assoc.}
\bvolume{118}
\bpages{1923--1939}.
\end{barticle}
\endbibitem

\bibitem[\protect\citeauthoryear{Hallin, La~Vecchia and
  Liu}{2023}]{HallinBernou23}
\begin{barticle}[author]
\bauthor{\bsnm{Hallin},~\bfnm{M.}\binits{M.}},
  \bauthor{\bsnm{La~Vecchia},~\bfnm{D.}\binits{D.}} \AND
  \bauthor{\bsnm{Liu},~\bfnm{H.}\binits{H.}}
(\byear{2023}).
\btitle{Rank-based testing for semiparametric VAR models: a measure
  transportation approach}.
\bjournal{Bernoulli}
\bvolume{29}
\bpages{229--273}.
\end{barticle}
\endbibitem

\bibitem[\protect\citeauthoryear{Hallin, Liu and
  Verdebout}{2024}]{VerdeboutJRSSB}
\begin{barticle}[author]
\bauthor{\bsnm{Hallin},~\bfnm{M.}\binits{M.}},
  \bauthor{\bsnm{Liu},~\bfnm{H.}\binits{H.}} \AND
  \bauthor{\bsnm{Verdebout},~\bfnm{T.}\binits{T.}}
(\byear{2024}).
\btitle{Nonparametric measure-transportation-based methods for directional
  data}.
\bjournal{J. R. Stat. Soc. Ser. B. Stat. Methodol.}
\bvolume{\textcolor{black}{86}}
\bpages{\textcolor{black}{1172--1196}}.
\end{barticle}
\endbibitem

\bibitem[\protect\citeauthoryear{Hallin et~al.}{2021}]{Hal21}
\begin{barticle}[author]
\bauthor{\bsnm{Hallin},~\bfnm{M.}\binits{M.}},
  \bauthor{\bsnm{Del~Barrio},~\bfnm{E.}\binits{E.}},
  \bauthor{\bsnm{J.},~\bfnm{Cuesta-Albertos}\binits{C.-A.}} \AND
  \bauthor{\bsnm{Matr\'{a}n},~\bfnm{C.}\binits{C.}}
(\byear{2021}).
\btitle{Distribution and quantile functions, ranks and signs in dimension $d$:
  A measure transportation approach}.
\bjournal{Ann. Statist.}
\bvolume{49}
\bpages{1139--1165}.
\end{barticle}
\endbibitem

\bibitem[\protect\citeauthoryear{Hampel}{1968}]{Hampel68}
\begin{bphdthesis}[author]
\bauthor{\bsnm{Hampel},~\bfnm{Frank~R.}\binits{F.~R.}}
(\byear{1968}).
\btitle{Contributions to the theory of robust estimation},
\btype{PhD thesis},
\bpublisher{Dept. Statistics, Univ. California, Berkeley}.
\end{bphdthesis}
\endbibitem

\bibitem[\protect\citeauthoryear{Hampel}{1971}]{Hampel71}
\begin{barticle}[author]
\bauthor{\bsnm{Hampel},~\bfnm{Frank~R.}\binits{F.~R.}}
(\byear{1971}).
\btitle{A general qualitative definition of robustness}.
\bjournal{Ann. Statist.}
\bvolume{42}
\bpages{1887--1896}.
\end{barticle}
\endbibitem

\bibitem[\protect\citeauthoryear{Hampel et~al.}{1986}]{Hametal1986}
\begin{bbook}[author]
\bauthor{\bsnm{Hampel},~\bfnm{Frank~R.}\binits{F.~R.}},
  \bauthor{\bsnm{Ronchetti},~\bfnm{Elvezio~M.}\binits{E.~M.}},
  \bauthor{\bsnm{Rousseeuw},~\bfnm{Peter~J.}\binits{P.~J.}} \AND
  \bauthor{\bsnm{Stahel},~\bfnm{Werner~A.}\binits{W.~A.}}
(\byear{1986}).
\btitle{Robust Statistics: The Approach Based on Influence Functions}.
\bseries{Wiley Series in Probability and Mathematical Statistics: Probability
  and Mathematical Statistics}.
\bpublisher{John Wiley \& Sons Inc.}, \baddress{New York}.
\end{bbook}
\endbibitem

\bibitem[\protect\citeauthoryear{Hundrieser, Klatt and
  Munk}{2024}]{Hundrieser2024}
\begin{barticle}[author]
\bauthor{\bsnm{Hundrieser},~\bfnm{S.}\binits{S.}},
  \bauthor{\bsnm{Klatt},~\bfnm{M.}\binits{M.}} \AND
  \bauthor{\bsnm{Munk},~\bfnm{A.}\binits{A.}}
(\byear{2024}).
\btitle{Limit distributions and sensitivity analysis for empirical entropic
  optimal transport on countable spaces}.
\bjournal{Ann. Appl. Probab.}
\bvolume{34}
\bpages{1403--1468}.
\end{barticle}
\endbibitem

\bibitem[\protect\citeauthoryear{Hundrieser, Staudt and Munk}{2024}]{HSM24}
\begin{barticle}[author]
\bauthor{\bsnm{Hundrieser},~\bfnm{S.}\binits{S.}},
  \bauthor{\bsnm{Staudt},~\bfnm{T.}\binits{T.}} \AND
  \bauthor{\bsnm{Munk},~\bfnm{A.}\binits{A.}}
(\byear{2024}).
\btitle{Empirical optimal transport between different measures adapts to lower
  complexity}.
\bjournal{Ann. Inst. Henri Poincar\'{e} Probab. Stat.}
\bvolume{60}
\bpages{824--846}.
\end{barticle}
\endbibitem

\bibitem[\protect\citeauthoryear{Hundrieser et~al.}{2024}]{Hundrieser2024b}
\begin{barticle}[author]
\bauthor{\bsnm{Hundrieser},~\bfnm{S.}\binits{S.}},
  \bauthor{\bsnm{Klatt},~\bfnm{M.}\binits{M.}},
  \bauthor{\bsnm{Staudt},~\bfnm{T.}\binits{T.}} \AND
  \bauthor{\bsnm{Munk},~\bfnm{A.}\binits{A.}}
(\byear{2024}).
\btitle{A unifying approach to distributional limits for empirical optimal
  transport}.
\bjournal{Bernoulli}
\bvolume{\textcolor{black}{30}}
\bpages{\textcolor{black}{2846--2877}}.
\end{barticle}
\endbibitem

\bibitem[\protect\citeauthoryear{H\"{u}tter and Rigollet}{2021}]{Rigollet2021}
\begin{barticle}[author]
\bauthor{\bsnm{H\"{u}tter},~\bfnm{J.~C.}\binits{J.~C.}} \AND
  \bauthor{\bsnm{Rigollet},~\bfnm{P.}\binits{P.}}
(\byear{2021}).
\btitle{Minimax estimation of smooth optimal transport maps}.
\bjournal{Ann. Statist.}
\bvolume{49}
\bpages{1166--1194}.
\end{barticle}
\endbibitem

\bibitem[\protect\citeauthoryear{Kallenberg}{2021}]{Kal}
\begin{bbook}[author]
\bauthor{\bsnm{Kallenberg},~\bfnm{O.}\binits{O.}}
(\byear{2021}).
\btitle{Foundations of Modern Probability}.
\bpublisher{Springer}, \baddress{Switzerland}.
\end{bbook}
\endbibitem

\bibitem[\protect\citeauthoryear{Kantorovich}{1942}]{Kanto1942}
\begin{barticle}[author]
\bauthor{\bsnm{Kantorovich},~\bfnm{L.}\binits{L.}}
(\byear{1942}).
\btitle{On the translocation of masses}.
\bjournal{Doklady Akademii Nauk URSS}
\bvolume{37}
\bpages{7--8}.
\end{barticle}
\endbibitem

\bibitem[\protect\citeauthoryear{Kolouri et~al.}{2017}]{Kolouri2017}
\begin{barticle}[author]
\bauthor{\bsnm{Kolouri},~\bfnm{S.}\binits{S.}},
  \bauthor{\bsnm{Park},~\bfnm{S.~R.}\binits{S.~R.}},
  \bauthor{\bsnm{Thorpe},~\bfnm{M.}\binits{M.}},
  \bauthor{\bsnm{Slep\v{c}ev},~\bfnm{D.}\binits{D.}} \AND
  \bauthor{\bsnm{Rohde},~\bfnm{G.~K.}\binits{G.~K.}}
(\byear{2017}).
\btitle{Optimal mass transport: Signal processing and machine-learning
  applications}.
\bjournal{IEEE Signal Process. Mag.}
\bvolume{34}
\bpages{43--59}.
\end{barticle}
\endbibitem

\bibitem[\protect\citeauthoryear{Koltchinski}{1997}]{Kol1997}
\begin{barticle}[author]
\bauthor{\bsnm{Koltchinski},~\bfnm{V.~I.}\binits{V.~I.}}
(\byear{1997}).
\btitle{{M-estimation}, convexity and quantiles}.
\bjournal{Ann. Statist.}
\bvolume{25}
\bpages{435--477}.
\end{barticle}
\endbibitem

\bibitem[\protect\citeauthoryear{Konen and Paindaveine}{2025}]{AIHP}
\begin{barticle}[author]
\bauthor{\bsnm{Konen},~\bfnm{D.}\binits{D.}} \AND
  \bauthor{\bsnm{Paindaveine},~\bfnm{D.}\binits{D.}}
(\byear{2025}).
\btitle{On the robustness of spatial quantiles}.
\bjournal{Ann. Inst. Henri Poincar\'{e} Probab. Stat.,}
\bpages{to appear}.
\end{barticle}
\endbibitem

\bibitem[\protect\citeauthoryear{Le et~al.}{2021}]{Le2021}
\begin{barticle}[author]
\bauthor{\bsnm{Le},~\bfnm{K.}\binits{K.}},
  \bauthor{\bsnm{Nguyen},~\bfnm{H.}\binits{H.}},
  \bauthor{\bsnm{Nguyen},~\bfnm{Q.~M.}\binits{Q.~M.}},
  \bauthor{\bsnm{Pham},~\bfnm{T.}\binits{T.}},
  \bauthor{\bsnm{Bui},~\bfnm{H.}\binits{H.}} \AND
  \bauthor{\bsnm{Ho},~\bfnm{N.}\binits{N.}}
(\byear{2021}).
\btitle{On robust optimal transport: Computational complexity and barycenter
  computation}.
\bjournal{Adv. Neural Inf. Process. Syst.}
\bvolume{34}.
\end{barticle}
\endbibitem

\bibitem[\protect\citeauthoryear{L\'{e}vy}{2015}]{Lev}
\begin{barticle}[author]
\bauthor{\bsnm{L\'{e}vy},~\bfnm{B.}\binits{B.}}
(\byear{2015}).
\btitle{A numerical algorithm for $L_2$ semi-discrete optimal transport in 3D}.
\bjournal{ESAIM: Math. Model. Numer. Anal.}
\bvolume{49}
\bpages{1693--1715}.
\end{barticle}
\endbibitem

\bibitem[\protect\citeauthoryear{Loeper}{2005}]{Loeper2005}
\begin{barticle}[author]
\bauthor{\bsnm{Loeper},~\bfnm{G.}\binits{G.}}
(\byear{2005}).
\btitle{On the regularity of the polar factorization for time dependent maps}.
\bjournal{Calc. Var. Partial Differential Equations}
\bvolume{22}
\bpages{343--374}.
\end{barticle}
\endbibitem

\bibitem[\protect\citeauthoryear{Lopuha\"a and Rousseeuw}{1991}]{LopRou1991}
\begin{barticle}[author]
\bauthor{\bsnm{Lopuha\"a},~\bfnm{Hendrik~P.}\binits{H.~P.}} \AND
  \bauthor{\bsnm{Rousseeuw},~\bfnm{Peter~J.}\binits{P.~J.}}
(\byear{1991}).
\btitle{Breakdown points of affine equivariant estimators of multivariate
  location and covariance matrices}.
\bjournal{Ann. Statist.}
\bvolume{19}
\bpages{229--248}.
\end{barticle}
\endbibitem

\bibitem[\protect\citeauthoryear{Manole and Niles-Weed}{2024}]{Manole2024a}
\begin{barticle}[author]
\bauthor{\bsnm{Manole},~\bfnm{T.}\binits{T.}} \AND
  \bauthor{\bsnm{Niles-Weed},~\bfnm{J.}\binits{J.}}
(\byear{2024}).
\btitle{Sharp convergence rates for empirical optimal transport with smooth
  costs}.
\bjournal{Ann. Appl. Probab.}
\bvolume{34}
\bpages{1108--1135}.
\end{barticle}
\endbibitem

\bibitem[\protect\citeauthoryear{Manole et~al.}{2024a}]{Manole2024}
\begin{barticle}[author]
\bauthor{\bsnm{Manole},~\bfnm{T.}\binits{T.}},
  \bauthor{\bsnm{Balakrishnan},~\bfnm{S.}\binits{S.}},
  \bauthor{\bsnm{Niles-Weed},~\bfnm{J.}\binits{J.}} \AND
  \bauthor{\bsnm{Wasserman},~\bfnm{L.}\binits{L.}}
(\byear{2024}a).
\btitle{Plugin estimation of smooth optimal transport maps}.
\bjournal{Ann. Statist.}
\bvolume{\textcolor{black}{52}}
\bpages{\textcolor{black}{966--998}}.
\end{barticle}
\endbibitem

\bibitem[\protect\citeauthoryear{Manole et~al.}{2024b}]{ManoleBNW2024}
\begin{barticle}[author]
\bauthor{\bsnm{Manole},~\bfnm{T.}\binits{T.}},
  \bauthor{\bsnm{Balakrishnan},~\bfnm{S.}\binits{S.}},
  \bauthor{\bsnm{Niles-Weed},~\bfnm{J.}\binits{J.}} \AND
  \bauthor{\bsnm{Wasserman},~\bfnm{L.}\binits{L.}}
(\byear{2024}b).
\btitle{Central limit theorems for smooth optimal transport maps}.
\bjournal{arXiv:2312.12407}.
\end{barticle}
\endbibitem

\bibitem[\protect\citeauthoryear{Mass{{\'e}}}{2009}]{Mas2009}
\begin{barticle}[author]
\bauthor{\bsnm{Mass{{\'e}}},~\bfnm{Jean-Claude}\binits{J.-C.}}
(\byear{2009}).
\btitle{Multivariate trimmed means based on the {T}ukey depth}.
\bjournal{J. Statist. Plann. Inference}
\bvolume{139}
\bpages{366--384}.
\end{barticle}
\endbibitem

\bibitem[\protect\citeauthoryear{McCann}{1995}]{McCann}
\begin{barticle}[author]
\bauthor{\bsnm{McCann},~\bfnm{R.~J.}\binits{R.~J.}}
(\byear{1995}).
\btitle{Existence and uniqueness of monotone measure-preserving maps}.
\bjournal{Duke Math. J.}
\bvolume{80}
\bpages{309--323}.
\end{barticle}
\endbibitem

\bibitem[\protect\citeauthoryear{Merigot}{2011}]{Mer}
\begin{barticle}[author]
\bauthor{\bsnm{Merigot},~\bfnm{Q.}\binits{Q.}}
(\byear{2011}).
\btitle{A multiscale approach to optimal transport}.
\bjournal{Comput. Graph. Forum}
\bvolume{30}
\bpages{1583--1592}.
\end{barticle}
\endbibitem

\bibitem[\protect\citeauthoryear{Monge}{1781}]{Monge1781}
\begin{barticle}[author]
\bauthor{\bsnm{Monge},~\bfnm{G.}\binits{G.}}
(\byear{1781}).
\btitle{M\'{e}moire sur la th\'{e}orie des d\'{e}blais et des remblais}.
\bjournal{Histoire de l'Acad\'{e}mie Royale des Sciences de Paris}
\bpages{666--704}.
\end{barticle}
\endbibitem

\bibitem[\protect\citeauthoryear{Mukherjee et~al.}{2021}]{Mukherjee2021}
\begin{barticle}[author]
\bauthor{\bsnm{Mukherjee},~\bfnm{D.}\binits{D.}},
  \bauthor{\bsnm{Guha},~\bfnm{A.}\binits{A.}},
  \bauthor{\bsnm{Solomon},~\bfnm{J.}\binits{J.}},
  \bauthor{\bsnm{Sun},~\bfnm{Y.}\binits{Y.}} \AND
  \bauthor{\bsnm{Yurochkin},~\bfnm{M}\binits{M.}}
(\byear{2021}).
\btitle{Outlier-robust optimal transport}.
\bjournal{ICML}.
\end{barticle}
\endbibitem

\bibitem[\protect\citeauthoryear{Nagy, Sch\"{u}tt and Werner}{2019}]{NagySS}
\begin{barticle}[author]
\bauthor{\bsnm{Nagy},~\bfnm{S.}\binits{S.}},
  \bauthor{\bsnm{Sch\"{u}tt},~\bfnm{C.}\binits{C.}} \AND
  \bauthor{\bsnm{Werner},~\bfnm{E.~M.}\binits{E.~M.}}
(\byear{2019}).
\btitle{Halfspace depth and floating body}.
\bjournal{Stat. Surv.}
\bvolume{13}
\bpages{52--118}.
\end{barticle}
\endbibitem

\bibitem[\protect\citeauthoryear{Nietert, Cummings and
  Goldfeld}{2022}]{Nietert2022}
\begin{barticle}[author]
\bauthor{\bsnm{Nietert},~\bfnm{S.}\binits{S.}},
  \bauthor{\bsnm{Cummings},~\bfnm{R.}\binits{R.}} \AND
  \bauthor{\bsnm{Goldfeld},~\bfnm{Z.}\binits{Z.}}
(\byear{2022}).
\btitle{Outlier-robust optimal transport: duality, structure, and statistical
  analysis}.
\bjournal{AISTATS}
\bvolume{151}.
\end{barticle}
\endbibitem

\bibitem[\protect\citeauthoryear{Orlova et~al.}{2016}]{Orlova16}
\begin{barticle}[author]
\bauthor{\bsnm{Orlova},~\bfnm{D.~Y.}\binits{D.~Y.}},
  \bauthor{\bsnm{Zimmerman},~\bfnm{N.}\binits{N.}},
  \bauthor{\bsnm{Meehan},~\bfnm{S.}\binits{S.}},
  \bauthor{\bsnm{Meehan},~\bfnm{C.}\binits{C.}},
  \bauthor{\bsnm{Waters},~\bfnm{J.}\binits{J.}},
  \bauthor{\bsnm{Ghosn},~\bfnm{E.~E.~B.}\binits{E.~E.~B.}},
  \bauthor{\bsnm{Filatenkov},~\bfnm{A.}\binits{A.}},
  \bauthor{\bsnm{Kolyagin},~\bfnm{G.~A.}\binits{G.~A.}},
  \bauthor{\bsnm{Gernez},~\bfnm{Y.}\binits{Y.}},
  \bauthor{\bsnm{Tsuda},~\bfnm{S.}\binits{S.}},
  \bauthor{\bsnm{Moore},~\bfnm{W.}\binits{W.}},
  \bauthor{\bsnm{Moss},~\bfnm{R.~B.}\binits{R.~B.}},
  \bauthor{\bsnm{Herzenberg},~\bfnm{L.~A.}\binits{L.~A.}} \AND
  \bauthor{\bsnm{Walther},~\bfnm{G.~A. Irving~C.}\binits{G.~A. I.~C.}}
(\byear{2016}).
\btitle{Earth Mover's Distance (EMD): A true metric for comparing biomarker
  expression levels in cell populations}.
\bjournal{PLoS ONE}
\bvolume{11}.
\end{barticle}
\endbibitem

\bibitem[\protect\citeauthoryear{Panaretos and Zemel}{2020}]{Panaretos2020}
\begin{bbook}[author]
\bauthor{\bsnm{Panaretos},~\bfnm{V.~M.}\binits{V.~M.}} \AND
  \bauthor{\bsnm{Zemel},~\bfnm{Y.}\binits{Y.}}
(\byear{2020}).
\btitle{An Invitation to Statistics in Wasserstein Space}.
\bpublisher{Springer Nature}.
\end{bbook}
\endbibitem

\bibitem[\protect\citeauthoryear{Peyr\'{e} and Cuturi}{2019}]{PeyrCutu2019}
\begin{barticle}[author]
\bauthor{\bsnm{Peyr\'{e}},~\bfnm{G.}\binits{G.}} \AND
  \bauthor{\bsnm{Cuturi},~\bfnm{M.}\binits{M.}}
(\byear{2019}).
\btitle{Computational optimal transport: With applications to data science}.
\bjournal{Found. Trends Mach. Learn.}
\bvolume{11}
\bpages{355--607}.
\end{barticle}
\endbibitem

\bibitem[\protect\citeauthoryear{Rockafellar}{1970}]{Rock1970}
\begin{bbook}[author]
\bauthor{\bsnm{Rockafellar},~\bfnm{R.~T.}\binits{R.~T.}}
(\byear{1970}).
\btitle{Convex Analysis}.
\bpublisher{Princeton University Press}.
\end{bbook}
\endbibitem

\bibitem[\protect\citeauthoryear{Ronchetti}{2023}]{Ronch2023}
\begin{binproceedings}[author]
\bauthor{\bsnm{Ronchetti},~\bfnm{Elvezio~M.}\binits{E.~M.}}
(\byear{2023}).
\btitle{Robustness aspects of optimal transport}.
In \bbooktitle{Research Papers in Statistical Inference for Time Series and
  Related Models. Essays in Honor of Masanobu Taniguchi}
(\beditor{\bfnm{Yan}\binits{Y.}~\bsnm{Liu}},
  \beditor{\bfnm{Junichi}\binits{J.}~\bsnm{Hirukawa}} \AND
  \beditor{\bfnm{Yoshihide}\binits{Y.}~\bsnm{Kakizawa}}, eds.).
\bpublisher{Springer Nature}, \baddress{Singapore}.
\end{binproceedings}
\endbibitem

\bibitem[\protect\citeauthoryear{Rousseeuw and Ruts}{1999}]{RouRut1999}
\begin{barticle}[author]
\bauthor{\bsnm{Rousseeuw},~\bfnm{Peter~J.}\binits{P.~J.}} \AND
  \bauthor{\bsnm{Ruts},~\bfnm{Ida}\binits{I.}}
(\byear{1999}).
\btitle{The depth function of a population distribution}.
\bjournal{Metrika}
\bvolume{49}
\bpages{213--244}.
\end{barticle}
\endbibitem

\bibitem[\protect\citeauthoryear{Rousseeuw and Struyf}{2004}]{RouStr2004}
\begin{barticle}[author]
\bauthor{\bsnm{Rousseeuw},~\bfnm{Peter~J.}\binits{P.~J.}} \AND
  \bauthor{\bsnm{Struyf},~\bfnm{Anja}\binits{A.}}
(\byear{2004}).
\btitle{Characterizing angular symmetry and regression symmetry}.
\bjournal{J. Statist. Plann. Inference}
\bvolume{122}
\bpages{161--173}.
\end{barticle}
\endbibitem

\bibitem[\protect\citeauthoryear{Sadhu, Goldfeld and Kato}{2024}]{Sadhu2023}
\begin{barticle}[author]
\bauthor{\bsnm{Sadhu},~\bfnm{R.}\binits{R.}},
  \bauthor{\bsnm{Goldfeld},~\bfnm{Z.}\binits{Z.}} \AND
  \bauthor{\bsnm{Kato},~\bfnm{K.}\binits{K.}}
(\byear{2024}).
\btitle{Stability and statistical inference for semidiscrete optimal transport
  maps}.
\bjournal{\textcolor{black}{Ann. Appl. Probab.}}
\bvolume{\textcolor{black}{34}}
\bpages{\textcolor{black}{5694--5736}}.
\end{barticle}
\endbibitem

\bibitem[\protect\citeauthoryear{Santambrogio}{2015}]{Santa2015}
\begin{bbook}[author]
\bauthor{\bsnm{Santambrogio},~\bfnm{F.}\binits{F.}}
(\byear{2015}).
\btitle{Optimal Transport for Applied Mathematicians: Calculus of Variations,
  PDEs, and Modeling}.
\bseries{Progress in Nonlinear Differential Equations and Their Applications}.
\bpublisher{Springer}, \baddress{Cham, Switzerland}.
\end{bbook}
\endbibitem

\bibitem[\protect\citeauthoryear{Sch\"{u}tt}{1991}]{Schu1991}
\begin{barticle}[author]
\bauthor{\bsnm{Sch\"{u}tt},~\bfnm{C.}\binits{C.}}
(\byear{1991}).
\btitle{The convex floating body and polyhedral approximation}.
\bjournal{Israel J. Math.}
\bvolume{73}
\bpages{65--77}.
\end{barticle}
\endbibitem

\bibitem[\protect\citeauthoryear{Segers}{2022}]{Segers2022}
\begin{barticle}[author]
\bauthor{\bsnm{Segers},~\bfnm{J.}\binits{J.}}
(\byear{2022}).
\btitle{Graphical and uniform consistency of estimated optimal transport
  plans}.
\bjournal{arXiv:2208.02508}.
\end{barticle}
\endbibitem

\bibitem[\protect\citeauthoryear{Shi, Drton and Han}{2022}]{Shi2022}
\begin{barticle}[author]
\bauthor{\bsnm{Shi},~\bfnm{H.}\binits{H.}},
  \bauthor{\bsnm{Drton},~\bfnm{M.}\binits{M.}} \AND
  \bauthor{\bsnm{Han},~\bfnm{F.}\binits{F.}}
(\byear{2022}).
\btitle{Distribution-free consistent independence tests via center-outward
  ranks and signs}.
\bjournal{J. Amer. Statist. Assoc.}
\bvolume{117}
\bpages{395--410}.
\end{barticle}
\endbibitem

\bibitem[\protect\citeauthoryear{Shi et~al.}{2022}]{ShiAoS22}
\begin{barticle}[author]
\bauthor{\bsnm{Shi},~\bfnm{H.}\binits{H.}},
  \bauthor{\bsnm{Hallin},~\bfnm{M.}\binits{M.}},
  \bauthor{\bsnm{Drton},~\bfnm{M.}\binits{M.}} \AND
  \bauthor{\bsnm{Han},~\bfnm{F.}\binits{F.}}
(\byear{2022}).
\btitle{On universally consistent and fully distribution-free rank tests of
  vector independence}.
\bjournal{Ann. Statist.}
\bvolume{50}
\bpages{1933--1959}.
\end{barticle}
\endbibitem

\bibitem[\protect\citeauthoryear{Shi et~al.}{2025}]{Shi24}
\begin{barticle}[author]
\bauthor{\bsnm{Shi},~\bfnm{H.}\binits{H.}},
  \bauthor{\bsnm{Drton},~\bfnm{M.}\binits{M.}},
  \bauthor{\bsnm{Hallin},~\bfnm{M.}\binits{M.}} \AND
  \bauthor{\bsnm{Han},~\bfnm{F.}\binits{F.}}
(\byear{2025}).
\btitle{Distribution-free tests of multivariate independence based on
  center-outward quadrant, Spearman, Kendall, and van der Waerden statistics}.
\bjournal{Bernoulli}
\bvolume{\textcolor{black}{31}}
\bpages{\textcolor{black}{106--129}}.
\end{barticle}
\endbibitem

\bibitem[\protect\citeauthoryear{Staudt and Hundrieser}{2025}]{Staudt2023}
\begin{barticle}[author]
\bauthor{\bsnm{Staudt},~\bfnm{T.}\binits{T.}} \AND
  \bauthor{\bsnm{Hundrieser},~\bfnm{S.}\binits{S.}}
(\byear{2025}).
\btitle{Convergence of empirical optimal transport in unbounded settings}.
\bjournal{Bernoulli}
\bvolume{\textcolor{black}{31}}
\bpages{\textcolor{black}{1929--1954}}.
\end{barticle}
\endbibitem

\bibitem[\protect\citeauthoryear{Tukey}{1975}]{Tuk1975}
\begin{binproceedings}[author]
\bauthor{\bsnm{Tukey},~\bfnm{John~W.}\binits{J.~W.}}
(\byear{1975}).
\btitle{Mathematics and the picturing of data}.
In \bbooktitle{Proceedings of the {I}nternational {C}ongress of
  {M}athematicians ({V}ancouver, {B}. {C}., 1974), {V}ol. 2}
\bpages{523--531}.
\bpublisher{Canad. Math. Congress, Montreal, Que.}
\end{binproceedings}
\endbibitem

\bibitem[\protect\citeauthoryear{Villani}{2008}]{Villani2008}
\begin{bbook}[author]
\bauthor{\bsnm{Villani},~\bfnm{C.}\binits{C.}}
(\byear{2008}).
\btitle{Optimal Transport: Old and New}.
\bseries{A Series of Comprehensive Studies in Mathematics}.
\bpublisher{Springer}, \baddress{Heidelberg}.
\end{bbook}
\endbibitem

\bibitem[\protect\citeauthoryear{Weed and Bach}{2019}]{WeedBach2019}
\begin{barticle}[author]
\bauthor{\bsnm{Weed},~\bfnm{J.}\binits{J.}} \AND
  \bauthor{\bsnm{Bach},~\bfnm{F.}\binits{F.}}
(\byear{2019}).
\btitle{Sharp asymptotic and finite-sample rates of convergence of empirical
  measures in Wasserstein distance}.
\bjournal{Bernoulli}
\bvolume{25}
\bpages{2620--2648}.
\end{barticle}
\endbibitem

\bibitem[\protect\citeauthoryear{Yang and Wang}{2024}]{Tengyao2025}
\begin{barticle}[author]
\bauthor{\bsnm{Yang},~\bfnm{X.}\binits{X.}} \AND
  \bauthor{\bsnm{Wang},~\bfnm{T.}\binits{T.}}
(\byear{2024}).
\btitle{\textcolor{black}{Multiple-output composite quantile regression through
  an optimal transport lens}}.
\bjournal{\textcolor{black}{PMLR}}
\bvolume{\textcolor{black}{247}}
\bpages{\textcolor{black}{5076--5122}}.
\end{barticle}
\endbibitem

\bibitem[\protect\citeauthoryear{Zuo}{2006}]{Zuo2006}
\begin{barticle}[author]
\bauthor{\bsnm{Zuo},~\bfnm{Yijun}\binits{Y.}}
(\byear{2006}).
\btitle{Multidimensional trimming based on projection depth}.
\bjournal{Ann. Statist.}
\bvolume{34}
\bpages{2211--2251}.
\end{barticle}
\endbibitem

\bibitem[\protect\citeauthoryear{Zuo and Serfling}{2000a}]{ZuoSer2000D}
\begin{barticle}[author]
\bauthor{\bsnm{Zuo},~\bfnm{Yijun}\binits{Y.}} \AND
  \bauthor{\bsnm{Serfling},~\bfnm{Robert}\binits{R.}}
(\byear{2000}a).
\btitle{On the performance of some robust nonparametric location measures
  relative to a general notion of multivariate symmetry}.
\bjournal{J. Statist. Plann. Inference}
\bvolume{84}
\bpages{55--79}.
\end{barticle}
\endbibitem

\bibitem[\protect\citeauthoryear{Zuo and Serfling}{2000b}]{ZuoSer2000A}
\begin{barticle}[author]
\bauthor{\bsnm{Zuo},~\bfnm{Yijun}\binits{Y.}} \AND
  \bauthor{\bsnm{Serfling},~\bfnm{Robert}\binits{R.}}
(\byear{2000}b).
\btitle{General notions of statistical depth function}.
\bjournal{Ann. Statist.}
\bvolume{28}
\bpages{461--482}.
\end{barticle}
\endbibitem

\end{thebibliography}

%
%
%

\end{document}